# NUMERICAL APPROXIMATIONS FOR A PHASE-FIELD MOVING CONTACT LINE MODEL WITH VARIABLE DENSITIES AND VISCOSITIES

HAIJUN YU† AND XIAOFENG YANG‡⋆

ABSTRACT. We consider the numerical approximations of a two-phase hydrodynamics coupled phase-field model that incorporates the variable densities, viscosities and moving contact line boundary conditions. The model is a nonlinear, coupled system that consists of incompressible Navier–Stokes equations with the generalized Navier boundary condition, and the Cahn–Hilliard equations with moving contact line boundary conditions. By some subtle explicit–implicit treatments to nonlinear terms, we develop two efficient, unconditionally energy stable numerical schemes, in particular, a linear decoupled energy stable scheme for the system with static contact line condition, and a nonlinear energy stable scheme for the system with dynamic contact line condition. An efficient spectral-Galerkin spatial discretization is implemented to verify the accuracy and efficiency of proposed schemes. Various numerical results show that the proposed schemes are efficient and accurate.

## 1. INTRODUCTION

Mixtures of two or more immiscible fluid components with different physical properties are widely used in many science and engineering applications. When the interface of a fluid mixture touches the solid wall, a physical process called "moving contact line" (MCL) occurs. Appearing in many applications (e.g., spray cooling of surfaces, crop spraying, spray coating, etc.), MCL problem has always been an appealing and challenging topic for mathematical modeling and simulations. Different to hydrodynamics of one simple fluid, the no-slip boundary condition for Navier–Stokes equations is not applicable for multi-fluid MCL problems because a non-physical velocity discontinuity will occur at the MCLs (cf. [10, 11, 36]). To understand the hydrodynamical behavior near the MCLs, several methods have been developed including MD simulations [29, 30, 56], microscopic-macroscopic hybrid simulations [21, 42], the level set method [55, 59, 60], the VOF method [26, 43], the front tracking method [24], and the phase-field approach [13, 27, 37, 38, 40, 44, 54] considered in this paper.

Among aforementioned models/numerical methods, the phase-field (or diffuse interface) approach is now popular and widely used to simulate the interfacial dynamics due to its versatility in modeling as well as simulations (cf. [4, 20, 22, 22, 25, 28, 31, 33, 34, 62, 65–67, 69], and the references therein). Its idea can actually be dated back to the ancient work of Rayleigh [41] and van der Waals [57] one century ago. Such a method considers the fluid-fluid interface as a continuous, but steep change of some physical properties of two fluids, e.g., density or viscosity, etc. An order parameter (or called phase field variable) is introduced to label the two fluid components, thus the fluid-fluid interface is then represented by a thin but smooth transition layer that can remove the singularities in practice. The model is then derived from the energy-based variational formalisms, thus the interfacial dynamics and the complex rheology are incorporated to an unified theoretical framework that allows the dynamical model for each component to be combined in a single model. Therefore, the developed governing system is normally well-posed and satisfies a thermodynamically consistent energy dissipation law (or called energy stable), that makes it possible to implement corresponding mathematical analysis or design efficient numerical methods.

We recall that a series of pioneering works about the macroscopic phase field modeling for the MCL problem, as well as their analysis and numerical simulations had been carried out by Qian et al. in [38–40]. The governing system consists of the Navier–Stokes equations with the general Navier boundary condition (GNBC), and the Cahn–Hilliard equation with the dynamic contact line condition (DCLC). From the numerical point of view, it is quite a challenging topic to develop efficient time marching schemes, in particular, the energy stable schemes to solve such a complex dynamical system. The difficulties include (i) the small interfacial width that introduces stiffness into the system; (ii) the nonlinear coupling between the phase variable and the velocity due to nonlinear convections and stresses; (iii) the nonlinear couplings between the velocity and the phase variable on the MCL boundaries; and (iv) the coupling between the density, the viscosity, the velocity and the pressure in the fluid momentum equation. Recently, several attempts were made to improve the numerical stability and efficiency of schemes for solving MCL problems including the work of He et al. [23], Gao and Wang [13, 14], Salgado et al. [44], Aland and Chen [3], Dong [8], Dong and Shen [9], and Shen et al. [54], etc. However, for the variable density and/or viscosity case, those schemes are either nonlinearly coupled [13, 14, 44] that requires some efficient iterative solvers and need relatively high computational cost, or linearly decoupled but unable to provide the energy stability in theory [8, 9].

Therefore, in this paper, we aim to construct some effective and efficient numerical schemes to solve the phase-field MCL model for the case of variable densities and viscosities. More precisely, the schemes are expected to, be unconditionally energy stable, satisfy an energy law in the discrete level, and lead to linear, decoupled, or coupled elliptic equations to solve at each time step.

The rest of the paper is organized as follows. In Section 2, we introduce the hydrodynamics coupled phase-field model with MCLs in the presence of non-matched densities and viscosities, and derive its associated PDE energy dissipation law. In Section 3, we present the numerical schemes, and prove their discrete energy dissipation law rigorously. In Section 4, we present the spatial discretization method using the Galerkin approach. In section 5, we present various numerical examples to illustrate the accuracy and efficiency of the proposed schemes. Some concluding remarks are given in Section 6.

## 2. The PDE system and its energy law.

We now describe the phase-field model for a mixture of two immiscible, incompressible fluids in a confined domain $\Omega \subset \mathbb{R}^d$ ($d = 2, 3$) with densities $\rho_1$, $\rho_2$ and viscosities $\mu_1$, $\mu_2$, respectively.

We introduce a phase field variable (macroscopic labeling function) $\phi(\boldsymbol{x}, t)$ such that

$$(2.1) \qquad \phi(\boldsymbol{x}, t) = \begin{cases} 1, & \text{fluid I}, \\ -1, & \text{fluid II}, \end{cases}$$

with a thin, smooth transition region of width $O(\epsilon)$, and consider the following Ginzburg-Landau type of Helmholtz free energy functional:

$$(2.2) \qquad E_{mix} = \lambda \int_\Omega \left(\frac{\epsilon}{2}|\nabla \phi|^2 + F(\phi)\right) d\bm{x},$$

where $\lambda$ denotes rescaled characteristic strength of phase mixing energy. The first term in $E_{mix}$ contributes to the hydro-philic type (tendency of mixing) of interactions between the materials and the second part, the double well bulk energy $F(\phi) = \frac{1}{4\epsilon}(\phi^2 - 1)^2$, represents the hydro-phobic type (tendency of separation) of interactions. As the consequence of the competition between the two types of interactions, the equilibrium configuration will include a diffusive interface with a thickness proportional to the parameter $\epsilon$ (cf., for instance, [67]).

The total energy of the hydrodynamic system is a sum of the kinetic energy $E_k$ together with the mixing energy $E_{mix}$:

$$(2.3) \qquad E = E_k + E_{mix} = \int_\Omega \left(\frac{\rho}{2}|\bm{u}|^2 + \lambda\big(\frac{\epsilon}{2}|\nabla \phi|^2 + F(\phi)\big)\right) d\bm{x},$$

where $\bm{u}$ is the fluid velocity field and $\rho$ is the density.

The evolution of the phase function is governed by the Cahn–Hilliard phase equation in the conserved form.

$$(2.4) \qquad \phi_t + \nabla \cdot (\bm{u}\phi) = M\Delta\mu,$$
$$(2.5) \qquad \mu = \lambda(-\epsilon\Delta\phi + f(\phi)),$$

where $\mu$ is the chemical potential and $M$ is a mobility parameter related to the relaxation time scale, and $f(\phi) = F'(\phi) = \frac{1}{\epsilon}\phi(\phi^2 - 1)$.

The momentum equation (macroscopic force balance) for the hydrodynamics takes the usual form of the Navier-Stokes equation:

$$(2.6) \qquad \rho(\bm{u}_t + (\bm{u} \cdot \nabla)\bm{u}) = \nabla \cdot \sigma,$$

where the total stress $\sigma = \nu D(\bm{u}) - pI + \sigma_e$ with $D(\bm{u}) = \nabla\bm{u} + \nabla\bm{u}^T$ and $\sigma_e$ is the extra elastic stress induced by the microscopic internal energy.

If we assume the two fluids have matched density and viscosity, i.e., $\rho_1 = \rho_2 = 1, \nu_1 = \nu_2 = \nu$, then the fluid momentum equation can be simplified as follows.

$$(2.7) \qquad \begin{cases} \bm{u}_t + (\bm{u} \cdot \nabla)\bm{u} - \nu\Delta\bm{u} + \nabla p + \phi\nabla\mu = 0, \\ \nabla \cdot \bm{u} = 0. \end{cases}$$

Note that the stress term is written as $\phi\nabla\mu$ instead of the $-\mu\nabla\phi$ in most other references. This is due to $\phi\nabla\mu = \nabla(\phi\mu) - \mu\nabla\phi$, and the $\nabla(\phi\mu)$ can be absorbed into the pressure term. About the theoretical or numerical study on this phase-field model, we refer to [4, 20, 25, 31, 54] and the references therein.

We now consider the case that $\rho_1 \neq \rho_2$. If the density ratio is small ($\sim O(1)$), one could use the well-known Boussinesq approximation to model the effect of density difference by a gravitational force (cf. for instance [31]). When the density ratio is large such that the Boussinesq approximation is no longer valid, the situation becomes more complicated, and there exist several phase-field models derived from various considerations (cf. [1, 16, 19, 32, 33, 48, 50, 52, 53]).

Assuming the density function $\rho$ and viscosity $\nu$ has the following linear relations,

$$
\begin{cases}
\rho(\phi) = \dfrac{\rho_1 - \rho_2}{2}\phi + \dfrac{\rho_1 + \rho_2}{2}, \\
\nu(\phi) = \dfrac{\nu_1 - \nu_2}{2}\phi + \dfrac{\nu_1 + \nu_2}{2}.
\end{cases} \tag{2.8}
$$

Then the mass conservation equation $\rho_t + \nabla \cdot (\rho \boldsymbol{u}) \neq 0$. In other words, from (2.4), the density function $\rho(\phi)$ satisfy the following diffusive equation

$$\rho_t + \nabla \cdot (\rho \boldsymbol{u}) + \nabla \cdot J = 0, \tag{2.9}$$

where $J = -\dfrac{\rho_1 - \rho_2}{2} M \nabla \mu$. We can easily derive the following identites

$$
\begin{aligned}
\dfrac{d}{dt}(\rho, \dfrac{|\boldsymbol{u}|^2}{2}) &= (\rho \boldsymbol{u}_t, \boldsymbol{u}) + (\rho_t, \dfrac{|\boldsymbol{u}|^2}{2}) \\
&= (\rho \boldsymbol{u}_t, \boldsymbol{u}) - (\nabla \cdot (\rho \boldsymbol{u}) + \nabla \cdot J, \dfrac{|\boldsymbol{u}|^2}{2}) \\
&= (\rho \boldsymbol{u}_t + \rho \boldsymbol{u} \cdot \nabla \boldsymbol{u} + J \cdot \nabla \boldsymbol{u}, \boldsymbol{u}),
\end{aligned} \tag{2.10}
$$

where we emphasize that only the boundary condition of $\boldsymbol{u} \cdot \boldsymbol{n}|_\Gamma = \partial_{\boldsymbol{n}} \mu|_\Gamma = 0$ is needed where $\boldsymbol{n}$ is the outward normal on the domain boundary $\Gamma$. Throughout this paper, we assume the MCL boundary is on the straight part of the domain boundary $\Gamma$. Thus we derived the Cahn–Hilliard–Navier–Stokes (CHNS) system with non-matched density and viscosity as follows (cf. [1]).

$$\phi_t + \nabla \cdot (\boldsymbol{u}\phi) = M\Delta\mu, \tag{2.11a}$$

$$\mu = \lambda(-\epsilon\Delta\phi + f(\phi)), \tag{2.11b}$$

$$\rho(\boldsymbol{u}_t + (\boldsymbol{u} \cdot \nabla)\boldsymbol{u}) + J \cdot \nabla \boldsymbol{u} + \nabla p - \nabla \cdot \nu D(\boldsymbol{u}) + \phi \nabla \mu = 0, \tag{2.11c}$$

$$\nabla \cdot \boldsymbol{u} = 0, \tag{2.11d}$$

where $p$, $\rho$ and $\nu$ are the pressure, density and viscosity of the mixture, respectively.

If the fluid-fluid interface never touches the wall (domain boundary $\Gamma$), we can assume the easy boundary conditions that can erase all boundary terms.

$$\boldsymbol{u}|_\Gamma = 0, \quad \partial_{\boldsymbol{n}}\phi|_\Gamma = 0, \quad \partial_{\boldsymbol{n}}\mu|_\Gamma = 0. \tag{2.12}$$

When the fluid-fluid interface touches the wall, the MCLs problem appears, we then have the GNBC boundary conditions for the velocity (cf. [38, 40]),

$$\boldsymbol{u} \cdot \boldsymbol{n} = 0, \quad \text{on } \Gamma, \tag{2.13}$$

$$\nu\ell(\phi)(\boldsymbol{u}_\tau - \boldsymbol{u}_w) + \nu\partial_{\boldsymbol{n}}\boldsymbol{u}_\tau - \lambda L(\phi)\nabla_\tau\phi = 0, \quad \text{on } \Gamma, \tag{2.14}$$

and together with the DCLC boundary conditions for the phase field variable,

$$\partial_{\boldsymbol{n}}\mu = 0, \quad \text{on } \Gamma, \tag{2.15}$$

$$\phi_t + \boldsymbol{u}_\tau \cdot \nabla_\tau \phi = -\gamma L(\phi), \quad \text{on } \Gamma. \tag{2.16}$$

Here $\ell(\phi) \geq 0$ is a given coefficient function that is the ratio of domain length to the slip length, $\gamma$ is a boundary relaxation coefficient, $\boldsymbol{u}_w$ is the wall velocity, $\boldsymbol{u}_\tau$ is the tangential velocity along the boundary tangential direction $\tau$, $\nabla_\tau = \nabla - (\boldsymbol{n} \cdot \nabla)\boldsymbol{n}$ is the gradient along $\tau$. The function $L(\phi)$ is defined as

$$L(\phi) = \epsilon\partial_{\boldsymbol{n}}\phi + g'(\phi), \tag{2.17}$$

where $g(\phi)$ is the interfacial energy with

$$(2.18) \qquad g(\phi) = -\frac{\sqrt{2}}{3}\cos\theta_s \sin\left(\frac{\pi}{2}\phi\right),$$

and $\theta_s$ is the static contact angle. From (2.13), we have $\boldsymbol{u} = \boldsymbol{u}_\tau$ on boundary $\Gamma$.

About the validity of the phase field model with GNBC (2.14) and the DCLC (2.16), we refer to a series of work of Qian and Wang in [38–40]. When $\gamma \to +\infty$, the dynamic contact line condition (2.15)-(2.16) reduces to the static contact line condition (SCLC) as

$$(2.19) \qquad \partial_{\boldsymbol{n}}\mu = 0, \quad L(\phi) = 0, \quad \text{on } \Gamma,$$

and the GNBC (2.13)-(2.14) reduces to the Navier boundary condition (NBC) as

$$(2.20) \qquad \boldsymbol{u}\cdot\boldsymbol{n} = 0, \quad \ell(\phi)(\boldsymbol{u}_\tau - \boldsymbol{u}_w) + \partial_{\boldsymbol{n}}\boldsymbol{u}_\tau = 0, \quad \text{on } \Gamma.$$

If we further set $g'(\phi) \equiv 0$ (i.e. $\theta_s = \frac{\pi}{2}$), the SCLC reduces to a phase-field model without contact line effect. If we take $\ell(\phi) \to +\infty$ in equation (2.20), namely, the slip length is zero, then the NBC reduces to the traditional no-slip boundary condition.

We now derive the PDE energy dissipation law for the system ((2.11), (2.13)-(2.16)). Here and after, for any function $f, g \in H^1(\Omega)$, we use $(f, g)$ to denote $\int_\Omega fg d\boldsymbol{x}$, $(f, g)_\Gamma$ to denote $\int_\Gamma fg ds$, and $\|f\|^2 = (f, f)$ and $\|f\|_\Gamma^2 = (f, f)_\Gamma$.

**Theorem 1.** *The system CHNS-GNBC-DCLC ((2.11), (2.13)-(2.16)) is a dissipative system satisfying the following energy dissipation law*

$$(2.21) \qquad \frac{\mathrm{d}}{\mathrm{d}t}E_{tot} = -\frac{1}{2}\|\sqrt{\nu}D(\boldsymbol{u})\|^2 - M\|\nabla\mu\|^2 - \lambda\gamma\|L(\phi)\|_\Gamma^2 - \|\sqrt{\nu\ell(\phi)}\boldsymbol{u}_s\|_\Gamma^2 - (\nu\ell(\phi)\boldsymbol{u}_s, \boldsymbol{u}_w)_\Gamma,$$

*where $\boldsymbol{u}_s = \boldsymbol{u} - \boldsymbol{u}_w$ is the velocity slip on boundary $\Gamma$, $E_{tot} = E_k[\rho, \boldsymbol{u}] + E_b[\phi] + E_s[\phi]$, and*

$$(2.22) \qquad E_k[\rho, \boldsymbol{u}] = \left(\rho, \frac{1}{2}|\boldsymbol{u}|^2\right), \quad E_b[\phi] = \lambda\frac{\epsilon}{2}\|\nabla\phi\|^2 + \lambda(F(\phi), 1), \quad E_s[\phi] = \lambda(g(\phi), 1)_\Gamma.$$

*Proof.* By taking the inner product of equation (2.11c) with $\boldsymbol{u}$, using the identity (2.10) and the incompressible condition (2.11d), we get

$$(2.23) \qquad \begin{aligned} \frac{\mathrm{d}}{\mathrm{d}t}\left(\rho, \frac{1}{2}|\boldsymbol{u}|^2\right) &= (\nabla\cdot\nu D(\boldsymbol{u}), \boldsymbol{u}) - (\phi\nabla\mu, \boldsymbol{u}) \\ &= -\frac{1}{2}\|\sqrt{\nu}D(\boldsymbol{u})\|^2 + (\nu\partial_{\boldsymbol{n}}\boldsymbol{u}_\tau, \boldsymbol{u}_\tau)_\Gamma - (\phi\nabla\mu, \boldsymbol{u}). \end{aligned}$$

By taking the inner product of equation (2.11a) with $\mu$, and using boundary conditions (2.13) and (2.15), we get

$$(2.24) \qquad (\phi_t, \mu) - (\boldsymbol{u}\phi, \nabla\mu) = -M\|\nabla\mu\|^2.$$

By taking the inner product of equation (2.11b) with $-\phi_t$, we have

$$(2.25) \qquad -(\mu, \phi_t) = \lambda\epsilon(\partial_{\boldsymbol{n}}\phi, \phi_t)_\Gamma - \frac{\lambda\epsilon}{2}\frac{\mathrm{d}}{\mathrm{d}t}\|\nabla\phi\|^2 - \lambda\frac{\mathrm{d}}{\mathrm{d}t}(F(\phi), 1).$$

Summing up equations (2.23)–(2.25), we obtain

$$(2.26) \qquad \begin{aligned} \frac{\mathrm{d}}{\mathrm{d}t}\left(\rho, \frac{1}{2}|\boldsymbol{u}|^2\right) + \lambda\frac{\mathrm{d}}{\mathrm{d}t}\frac{\epsilon}{2}\|\nabla\phi\|^2 + \lambda\frac{\mathrm{d}}{\mathrm{d}t}(F(\phi), 1) \\ = -\frac{1}{2}\|\sqrt{\nu}D(\boldsymbol{u})\|^2 - M\|\nabla\mu\|^2 + (\nu\partial_{\boldsymbol{n}}\boldsymbol{u}, \boldsymbol{u})_\Gamma + \lambda\epsilon(\partial_{\boldsymbol{n}}\phi, \phi_t)_\Gamma. \end{aligned}$$

Then, we use boundary condition (2.14) and (2.16)-(2.17), to derive

$$(2.27) \quad \begin{aligned}(\nu\partial_{\boldsymbol{n}}\boldsymbol{u},\boldsymbol{u})_\Gamma = (\nu\partial_{\boldsymbol{n}}\boldsymbol{u}_\tau,\boldsymbol{u}_\tau)_\Gamma &= \big(\lambda L(\phi)\nabla_\tau\phi - \nu\ell(\phi)(\boldsymbol{u}_\tau - \boldsymbol{u}_w),\boldsymbol{u}_\tau\big)_\Gamma \\ &= \lambda\left(L(\phi)\nabla_\tau\phi,\boldsymbol{u}_\tau\right)_\Gamma - (\nu\ell(\phi)\boldsymbol{u}_s,\boldsymbol{u}_s + \boldsymbol{u}_w)_\Gamma,\end{aligned}$$

$$(2.28) \quad \begin{aligned}\lambda\epsilon(\partial_{\boldsymbol{n}}\phi,\phi_t)_\Gamma &= \lambda\left(L(\phi) - g'(\phi),\phi_t\right)_\Gamma \\ &= \lambda\left(L(\phi),\phi_t\right)_\Gamma - \lambda\left(g'(\phi),\phi_t\right)_\Gamma \\ &= \lambda\left(L(\phi),-\boldsymbol{u}_\tau\cdot\nabla_\tau\phi - \gamma L(\phi)\right)_\Gamma - \lambda\frac{\mathrm{d}}{\mathrm{d}t}\left(g(\phi),1\right)_\Gamma \\ &= -\lambda\left(L(\phi)\nabla_\tau\phi,\boldsymbol{u}_\tau\right)_\Gamma - \lambda\gamma\|L(\phi)\|_\Gamma^2 - \lambda\frac{\mathrm{d}}{\mathrm{d}t}\left(g(\phi),1\right)_\Gamma.\end{aligned}$$

Summing up (2.26), (2.27) and (2.28), we get the desired energy law (2.21). □

**Remark 2.1.** *For the system of CHNS-NBC-SCLC ((2.11),(2.20) -(2.19)), we can derive the similar energy dissipative law.*

$$(2.29) \quad \frac{\mathrm{d}}{\mathrm{d}t}E_{tot} = -\frac{1}{2}\|\sqrt{\nu}D(\boldsymbol{u})\|^2 - M\|\nabla\mu\|^2 - \|\sqrt{\nu\ell(\phi)}\boldsymbol{u}_s\|_\Gamma^2 - (\nu\ell(\phi)\boldsymbol{u}_s,\boldsymbol{u}_w)_\Gamma.$$

## 3. Numerical Schemes.

We now present our numerical schemes to solve the coupling system. So far, there are two popular approaches to handle the Ginzburg–Landau potential $F(\phi)$. One is the convex splitting method (cf. [12, 70]), another is the stabilization method (cf. [49, 58, 72]). Here we adopt the latter one, that can provide a linear discretization for $f(\phi)$ instead of solving a nonlinear equation. The unconditional stability of the stabilization method requires that the second derivative of $F(\phi)$ to be bounded. However, this is not satisfied by the Ginzburg-Landau potential. Since we are only interested in $\phi \in [-1,1]$, and it is proved by [5] that a truncated $F(\phi)$ with quadratic growth at infinity also guarantees the boundless of $\phi$ in Cahn–Hilliard equation. So it is a common practice to modify $F(\phi)$ to have a quadratic growth rate for $|\phi| > 1$ (see e.g. [6, 7, 47–49, 61, 64, 68, 71]). Without loss of generality, we introduce the following $\widehat{F}(\phi)$ to replace $F(\phi)$:

$$(3.1) \quad \widehat{F}(\phi) = \frac{1}{4\epsilon}\begin{cases}4(\phi+1)^2, & \text{if } \phi < -1, \\ (\phi^2-1)^2, & \text{if } -1 \leq \phi \leq 1, \\ 4(\phi-1)^2, & \text{if } \phi > 1.\end{cases}$$

Correspondingly, we define $\widehat{f}(\phi) = \widehat{F}'(\phi)$ and

$$(3.2) \quad L_1 := \max_{\phi\in\mathbb{R}}|\widehat{f}'(\phi)| = \frac{2}{\epsilon}, \quad L_2 := \max_{\phi\in\mathbb{R}}|g''(\phi)| = \frac{\sqrt{2}\pi^2}{12}|\cos\theta_s|.$$

### 3.1. Linearly Decoupled Stable (LDS) Scheme.
We first study the system of CHNS-NBC-SCLC ((2.11),(2.20) -(2.19)) that is the simple version of the MCL phase-field model.

Define a cut-off function

$$(3.3) \quad \widehat{\phi} = \begin{cases}\phi, & |\phi| \leq 1, \\ \mathrm{sign}(\phi), & |\phi| > 1.\end{cases}$$

We use $\delta t$ to denote time step size, and a superscript $n$ on $\boldsymbol{u}$, $p$, $\phi$, $\mu$ to denote approximations of corresponding variables at time $n\delta t$.

Given $\rho^n$, $\nu^n$, $\boldsymbol{u}^n$, $\phi^n$, $p^n$, the LDS scheme calculate $\rho^{n+1}$, $\nu^{n+1}$, $\boldsymbol{u}^{n+1}$, $\phi^{n+1}$, $p^{n+1}$ and $\mu^{n+1}$ in four steps.

**Step 1:** Update $\phi^{n+1}$ and $\mu^{n+1}$ by solving

$$\frac{\phi^{n+1} - \phi^n}{\delta t} + \nabla \cdot (\boldsymbol{u}_*^n \phi^n) = M\Delta \mu^{n+1}, \tag{3.4}$$

$$\mu^{n+1} = \lambda\big(-\epsilon\Delta\phi^{n+1} + \widehat{f}(\phi^n) + S_1(\phi^{n+1} - \phi^n)\big), \tag{3.5}$$

with boundary conditions

$$\partial_{\boldsymbol{n}} \mu^{n+1} = 0, \quad \text{on } \Gamma, \tag{3.6}$$

$$\tilde{L}^{n+1} = 0, \quad \text{on } \Gamma, \tag{3.7}$$

where

$$\boldsymbol{u}_*^n = \boldsymbol{u}^n - \delta t \frac{\phi^n \nabla \mu^{n+1}}{\rho^n}, \tag{3.8}$$

$$\tilde{L}^{n+1} := \epsilon \partial_{\boldsymbol{n}} \phi^{n+1} + g'(\phi^n) + S_2(\phi^{n+1} - \phi^n), \tag{3.9}$$

where $S_1, S_2$ are two positive constants with the magnitude determined later.

**Step 2:** Update $\rho^{n+1}, \nu^{n+1}$ using

$$\rho^{n+1} = \frac{\rho_1 - \rho_2}{2} \widehat{\phi}^{n+1} + \frac{\rho_1 + \rho_2}{2}, \tag{3.10}$$

$$\nu^{n+1} = \frac{\nu_1 - \nu_2}{2} \widehat{\phi}^{n+1} + \frac{\nu_1 + \nu_2}{2}. \tag{3.11}$$

**Step 3:** Update $\boldsymbol{u}^{n+1}$ by solving

$$\begin{aligned}\rho^n \frac{\boldsymbol{u}^{n+1} - \boldsymbol{u}^n}{\delta t} - \nabla \cdot \nu^n D(\boldsymbol{u}^{n+1}) + \nabla(2p^n - p^{n-1}) + \rho^n \boldsymbol{u}^n \cdot \nabla \boldsymbol{u}^{n+1} + J^n \cdot \nabla \boldsymbol{u}^{n+1}, \\ + \phi^n \nabla \mu^{n+1} + \frac{1}{2}\frac{\rho^{n+1} - \rho^n}{\delta t} \boldsymbol{u}^{n+1} + \frac{1}{2}\nabla \cdot (\rho^n \boldsymbol{u}^n) \boldsymbol{u}^{n+1} + \frac{1}{2}\nabla \cdot J^n \boldsymbol{u}^{n+1} = 0,\end{aligned} \tag{3.12}$$

with boundary conditions

$$\boldsymbol{u}^{n+1} \cdot \boldsymbol{n} = 0, \quad \text{on } \Gamma, \tag{3.13}$$

$$\partial_{\boldsymbol{n}} \boldsymbol{u}_\tau^{n+1} + \ell(\phi^n) \boldsymbol{u}_s^{n+1} = 0, \quad \text{on } \Gamma, \tag{3.14}$$

where

$$J^n = -M \frac{\rho_1 - \rho_2}{2} \nabla \mu^n. \tag{3.15}$$

**Step 4:** Update $p^{n+1}$ by solving

$$\Delta(p^{n+1} - p^n) = \frac{\chi}{\delta t} \nabla \cdot \boldsymbol{u}^{n+1}, \tag{3.16}$$

with boundary conditions

$$\partial_{\boldsymbol{n}} p^{n+1} = 0 \quad \text{on } \Gamma, \tag{3.17}$$

and $\chi = \frac{1}{2}\min(\rho_1, \rho_2)$.

**Remark 3.1.** *Several remarks are in order.*

- *The last three terms in (3.12) is a first-order approximation of the term $\frac{1}{2}(\rho_t + \nabla \cdot (\rho\boldsymbol{u}) + \nabla \cdot J)\boldsymbol{u}$ at $t_{n+1}$. In the PDE system (2.11c), this term vanishes due to (2.9). Hence, (3.12) is indeed a consistent first-order approximation to (2.11c).*
- *We derive from (3.10), (3.11) and (3.3) that $\rho^{n+1} \geq \min(\rho_1, \rho_2)$ and $\mu^{n+1} \geq \min(\mu_1, \mu_2)$. In order to avoid solving an elliptic equation with $1/\rho$ as the variable coefficient, we adapt the pressure-stabilized scheme (3.16) to solve the pressure, which leads to a pressure Poisson equation (cf. [17,32,50]). There are several versions of projection type schemes to design the scheme for the pressure (cf. [15,18,46]). The pressure-stabilized formulation is specifically efficient to handle the variable density and/or viscosity problems.*
- *For stability, notice that the nonlinear functional $f(\phi)$ is proportional to $\frac{1}{\epsilon}$ with $\epsilon \ll 1$, thus the explicit treatment of this term usually leads to a restriction on the time step $\delta t$. It is a common practice to introduce a linear "stabilizing" term to improve the stability while preserving the simplicity (cf. [9,32,44,49–52,54,58,63,67]). It allows us to treat the nonlinear term explicitly without suffering from any time step constraint.*
- *For accuracy, we must notice that this stabilizing term introduces an extra consistent error of order $O(\frac{\delta t}{\epsilon})$ in a small region near the interface. This extra error is essentially the same order as the error induced by the explicit treatment of $f(\phi)$. (About the rigorous error analysis of the stabilized approach, we refer to [49]). Moreover, comparing to other nonlinear methods, e.g., convex splitting, the truncation error of the stabilized approach is essentially the same as that of the convex splitting [12], where the concave part of the nonlinear term ($\sim O(\frac{\delta t}{\epsilon})$) is treated explicitly as well. We must notice that for fixed time step, when $\epsilon$ becomes smaller, the extra stabilizing will definitely induce larger errors. Therefore, in practice, we usually fix the small parameter $\epsilon$ and then adjust the time step $\delta t$ to obtain the desired accuracy.*
- *We introduce an explicit velocities $\boldsymbol{u}_*^n$ in the convection term of phase field equation inspired by [35,51]. Such a first order approximation to velocity $\boldsymbol{u}^{n+1}$ in fact makes it possible to design a decoupled scheme while maintaining the energy stability.*
- *In the above scheme, computations of $(\phi^{n+1}, \mu^{n+1})$, $\boldsymbol{u}^{n+1}$, and $p^{n+1}$ are totally decoupled. Furthermore, each of the steps consists of solving a* linear *elliptic problem. To the best of authors' knowledge, this is the first linear, decoupled energy stable scheme for the phase-field model with MCLs for non-matched density and viscosity.*

**Theorem 2.** *Assuming $\boldsymbol{u}_w = 0$, and $S_1 \geq L_1/2$, $S_2 \geq L_2/2$, then the scheme (3.4)–(3.17) is energy stable in the sense that*

$$(3.18) \quad E^{n+1} \leq E^n - \frac{\delta t}{2}\|\sqrt{\nu^n}D(\boldsymbol{u}^{n+1})\|^2 - \delta t M\|\nabla\mu^{n+1}\|^2 - \delta t\|\sqrt{\nu^n\ell(\phi^n)}\boldsymbol{u}_s^{n+1}\|_\Gamma^2,$$

*where*

$$(3.19) \quad E^n = \frac{1}{2}\|\sqrt{\rho^n}\boldsymbol{u}^n\|^2 + \lambda\Big(\frac{\epsilon}{2}\|\nabla\phi^n\|^2 + (\widehat{F}(\phi^n), 1)\Big) + \frac{\delta t^2}{2\chi}\|\nabla p^n\|^2 + \lambda(g(\phi^n), 1)_\Gamma.$$

*Proof.* (i) Using intergration by parts and boundary condition $\boldsymbol{u} \cdot \boldsymbol{n}|_\Gamma = 0$, we can derive

$$(3.20) \quad (\boldsymbol{u} \cdot \nabla \boldsymbol{v}, \boldsymbol{v}) + \frac{1}{2}((\nabla \cdot \boldsymbol{u})\boldsymbol{v}, \boldsymbol{v}) = 0,$$

$$(3.21) \quad (\nabla p, \boldsymbol{u}) = -(\nabla \cdot \boldsymbol{u}, p).$$

Then we have

$$(3.22) \quad \left((\rho^n \boldsymbol{u}^n) \cdot \nabla \boldsymbol{u}^{n+1} + \frac{1}{2}\nabla \cdot (\rho^n \boldsymbol{u}^n)\boldsymbol{u}^{n+1}, \boldsymbol{u}^{n+1}\right) = 0,$$

$$(3.23) \quad \left(J^n \cdot \nabla \boldsymbol{u}^{n+1} + \frac{1}{2}(\nabla \cdot J^n)\boldsymbol{u}^{n+1}, \boldsymbol{u}^{n+1}\right) = 0.$$

By taking the $L^2$ inner product of (3.12) with $2\delta t \boldsymbol{u}^{n+1}$ and using an identity of

$$(3.24) \quad 2a(a-b) = a^2 - b^2 + (a-b)^2,$$

we can derive

$$(3.25) \quad \begin{aligned} &\|\sigma^n \boldsymbol{u}^{n+1}\|^2 - \|\sigma^n \boldsymbol{u}_*^n\|^2 + \|\sigma^n(\boldsymbol{u}^{n+1} - \boldsymbol{u}_*^n)\|^2 + \|\sigma^{n+1}\boldsymbol{u}^{n+1}\|^2 - \|\sigma^n \boldsymbol{u}^{n+1}\|^2 \\ &+ \delta t \|\sqrt{\nu^n} D(\boldsymbol{u}^{n+1})\|^2 - 2\delta t(\nu^n \partial_{\boldsymbol{n}}\boldsymbol{u}^{n+1}, \boldsymbol{u}^{n+1})_\Gamma \\ &+ 2\delta t(p^{n+1} - 2p^n + p^{n-1}, \nabla \cdot \boldsymbol{u}^{n+1}) - 2\delta t(p^{n+1}, \nabla \cdot \boldsymbol{u}^{n+1}) = 0, \end{aligned}$$

where $\sigma^n = \sqrt{\rho^n}$.

For the boundary term in (3.25), using equations (3.14), we derive

$$(3.26) \quad 2\delta t(\nu^n \partial_{\boldsymbol{n}}\boldsymbol{u}^{n+1}, \boldsymbol{u}^{n+1})_\Gamma = -2\delta t\|\sqrt{\nu^n \ell(\phi^n)}\boldsymbol{u}_s^{n+1}\|_\Gamma^2.$$

By taking the $L^2$ inner product of (3.16) with $\frac{2\delta t^2}{\chi}(p^{n+1} - 2p^n + p^{n-1})$ and with $-\frac{2\delta t^2}{\chi}p^{n+1}$ separately, we obtain

$$(3.27) \quad \begin{aligned} -\frac{\delta t^2}{\chi}\big(\|\nabla(p^{n+1} - p^n)\|^2 - \|\nabla(p^n - p^{n-1})\|^2 + \|\nabla(p^{n+1} - 2p^n + p^{n-1})\|^2\big) \\ = 2\delta t(\nabla \cdot \boldsymbol{u}^{n+1}, p^{n+1} - 2p^n + p^{n-1}), \end{aligned}$$

and

$$(3.28) \quad \frac{\delta t^2}{\chi}(\|\nabla p^{n+1}\|^2 - \|\nabla p^n\|^2 + \|\nabla(p^{n+1} - p^n)\|^2) = -2\delta t(\nabla \cdot \boldsymbol{u}^{n+1}, p^{n+1}).$$

Adding the above two equalities together, we get

$$(3.29) \quad \begin{aligned} &2\delta t(p^{n+1} - 2p^n + p^{n-1}, \nabla \cdot \boldsymbol{u}^{n+1}) - 2\delta t(p^{n+1}, \nabla \cdot \boldsymbol{u}^{n+1}) \\ &= \frac{\delta t^2}{\chi}(\|\nabla p^{n+1}\|^2 - \|\nabla p^n\|^2) + \frac{\delta t^2}{\chi}\|\nabla(p^n - p^{n-1})\|^2 - \frac{\delta t^2}{\chi}\|\nabla(p^{n+1} - 2p^n + p^{n-1})\|^2. \end{aligned}$$

Next, we take the difference of (3.16) at step $t^{n+1}$ and step $t^n$, pair the resulting equation with $p^{n+1} - 2p^n + p^{n-1}$, then take integration by parts for both sides to derive

$$(3.30) \quad \frac{\delta t^2}{\chi}\|\nabla(p^{n+1} - 2p^n + p^{n-1})\|^2 \leq \chi\|\boldsymbol{u}^{n+1} - \boldsymbol{u}^n\|^2 \leq \frac{1}{2}\|\sigma^n(\boldsymbol{u}^{n+1} - \boldsymbol{u}^n)\|^2.$$

Combining (3.25), (3.26), (3.29), and (3.30), we derive

$$(3.31) \quad \begin{aligned} &\|\sigma^{n+1}\boldsymbol{u}^{n+1}\|^2 - \|\sigma^n \boldsymbol{u}_*^n\|^2 + \|\sigma^n(\boldsymbol{u}^{n+1} - \boldsymbol{u}_*^n)\|^2 + \delta t\|\sqrt{\nu^n} D(\boldsymbol{u}^{n+1})\|^2 \\ &+ \frac{\delta t^2}{\chi}(\|\nabla p^{n+1}\|^2 - \|\nabla p^n\|^2) + \frac{\delta t^2}{\chi}\|\nabla(p^n - p^{n-1})\|^2 \\ &\leq \frac{1}{2}\|\sigma^n(\boldsymbol{u}^{n+1} - \boldsymbol{u}^n)\|^2 - 2\delta t\|\sqrt{\nu^n \ell(\phi^n)}\tilde{\boldsymbol{u}}_s^{n+1}\|_\Gamma^2. \end{aligned}$$

Noticing (3.8), we have

$$\frac{\rho^n(\boldsymbol{u}_*^n - \boldsymbol{u}^n)}{\delta t} = -\phi^n \nabla \mu^{n+1}, \tag{3.32}$$

By taking inner product of (3.32) with $2\delta t \boldsymbol{u}_*^n$, we have

$$\|\sigma^n \boldsymbol{u}_*^n\|^2 - \|\sigma^n \boldsymbol{u}^n\|^2 + \|\sigma^n(\boldsymbol{u}_*^n - \boldsymbol{u}^n)\|^2 = -2\delta t(\phi^n \nabla \mu^{n+1}, \boldsymbol{u}_*^n). \tag{3.33}$$

On the other hand, we derive from the triangle inequality that

$$\|\sigma^n(\boldsymbol{u}_*^n - \boldsymbol{u}^n)\|^2 + \|\sigma^n(\boldsymbol{u}^{n+1} - \boldsymbol{u}_*^n)\|^2 \geq \frac{1}{2}\|\sigma^n(\boldsymbol{u}^{n+1} - \boldsymbol{u}^n)\|^2. \tag{3.34}$$

Thus, by combining (3.31), (3.33), and (3.34), we obtain

$$\begin{aligned}
&\|\sigma^{n+1}\boldsymbol{u}^{n+1}\|^2 - \|\sigma^n \boldsymbol{u}^n\|^2 + \delta t\|\sqrt{\nu^n}D(\boldsymbol{u}^{n+1})\|^2 + \frac{\delta t^2}{\chi}(\|\nabla p^{n+1}\|^2 - \|\nabla p^n\|^2) \\
&+ \frac{\delta t^2}{\chi}\|\nabla(p^n - p^{n-1})\|^2 \leq -2\delta t(\phi^n \nabla \mu^{n+1}, \boldsymbol{u}_*^n) - 2\delta t\|\sqrt{\nu^n \ell(\phi^n)}\boldsymbol{u}_s^{n+1}\|_\Gamma^2.
\end{aligned} \tag{3.35}$$

(ii) Taking inner product of (3.4) with $\mu^{n+1}$, and using (3.13), (3.6), we have

$$2(\phi^{n+1} - \phi^n, \mu^{n+1}) - 2\delta t\left(\boldsymbol{u}_*^n \phi^n, \nabla \mu^{n+1}\right) = -2\delta t M \|\nabla \mu^{n+1}\|^2. \tag{3.36}$$

(iii) Taking inner product of (3.5) with $-2(\phi^{n+1} - \phi^n)$, we have

$$\begin{aligned}
-2(\mu^{n+1}, \phi^{n+1} - \phi^n) = &-\lambda\epsilon\left(\|\nabla\phi^{n+1}\|^2 - \|\nabla\phi^n\|^2 + \|\nabla(\phi^{n+1} - \phi^n)\|^2\right) \\
&+ 2\lambda\epsilon(\partial_{\boldsymbol{n}}\phi^{n+1}, \phi^{n+1} - \phi^n)_\Gamma \\
&- 2\lambda(\widehat{f}(\phi^n) + S_1(\phi^{n+1} - \phi^n), \phi^{n+1} - \phi^n).
\end{aligned} \tag{3.37}$$

For the boundary integral terms in (3.37), by using (3.7), we have

$$2\lambda\epsilon(\partial_{\boldsymbol{n}}\phi^{n+1}, \phi^{n+1} - \phi^n)_\Gamma = -2\lambda(g'(\phi^n) + S_2(\phi^{n+1} - \phi^n), \phi^{n+1} - \phi^n)_\Gamma. \tag{3.38}$$

By Taylor expansions of $\widehat{F}(\phi)$ and $g(\phi)$, we know there exist $\xi, \zeta$ such that

$$\widehat{f}(\phi^n)(\phi^{n+1} - \phi^n) = \widehat{F}(\phi^{n+1}) - \widehat{F}(\phi^n) - \frac{\widehat{f}'(\xi)}{2}(\phi^{n+1} - \phi^n)^2, \tag{3.39}$$

$$g'(\phi^n)(\phi^{n+1} - \phi^n) = g(\phi^{n+1}) - g(\phi^n) - \frac{g''(\zeta)}{2}(\phi^{n+1} - \phi^n)^2. \tag{3.40}$$

Combining equations (3.37), (3.38), (3.39) and (3.40), we get

$$\begin{aligned}
-2(\mu^{n+1}, \phi^{n+1} - \phi^n) = &-\lambda\epsilon(\|\nabla\phi^{n+1}\|^2 - \|\nabla\phi^n\|^2 + \|\nabla(\phi^{n+1} - \phi^n)\|^2) \\
&- 2\lambda(\widehat{F}(\phi^{n+1}) - \widehat{F}(\phi^n), 1) - 2\lambda\left(S_1 - \frac{\widehat{f}'(\xi)}{2}, (\phi^{n+1} - \phi^n)^2\right) \\
&- 2\lambda(g(\phi^{n+1}) - g(\phi^n), 1)_\Gamma - 2\lambda\left(S_2 - \frac{g''(\zeta)}{2}, (\phi^{n+1} - \phi^n)^2\right)_\Gamma.
\end{aligned} \tag{3.41}$$

(iv) Summing up equations (3.35), (3.36), and (3.41), we get

$$\frac{1}{2}\|\sigma^{n+1}\boldsymbol{u}^{n+1}\|^2 - \frac{1}{2}\|\sigma^n\boldsymbol{u}^n\|^2 + \lambda\left(\frac{\epsilon}{2}\|\nabla\phi^{n+1}\|^2 - \frac{\epsilon}{2}\|\nabla\phi^n\|^2 + \frac{\epsilon}{2}\|\nabla(\phi^{n+1} - \phi^n)\|^2\right)$$
$$+ \lambda(\widehat{F}(\phi^{n+1}) - \widehat{F}(\phi^n), 1) + \lambda(g(\phi^{n+1}) - g(\phi^n), 1)_\Gamma$$
$$+ \frac{\delta t^2}{2\chi}\left(\|\nabla p^{n+1}\|^2 - \|\nabla p^n\|^2\right) + \frac{\delta t^2}{2\chi}\|\nabla p^{n+1} - \nabla p^n\|^2$$
$$\leq -\frac{\delta t}{2}\|\sqrt{\nu^n}D(\boldsymbol{u}^{n+1})\|^2 - \delta t M\|\nabla\mu^{n+1}\|^2 - \delta t\|\sqrt{\nu^n\ell(\phi^n)}\boldsymbol{u}_s^{n+1}\|_\Gamma^2$$
$$- \lambda\left(S_1 - \frac{\widehat{f}'(\xi)}{2}, (\phi^{n+1} - \phi^n)^2\right) - \lambda\left(S_2 - \frac{g''(\zeta)}{2}, (\phi^{n+1} - \phi^n)^2\right)_\Gamma.$$

By the assumption $S_1 \geq L_1/2$, and $S_2 \geq L_2/2$, we get the desired energy estimate. □

3.2. **Nonlinear Coupled Scheme.** We now present a numerical scheme to solve the system of CHNS-GNBC-DCLC ((2.11), (2.13)-(2.16)).

Given $\rho^n, \nu^n, \boldsymbol{u}^n, \phi^n, p^n$, the scheme calculate $\rho^{n+1}, \nu^{n+1}, \boldsymbol{u}^{n+1}, \phi^{n+1}, p^{n+1}$ and $\mu^{n+1}$ in two steps.

*Step 1:* Update $\phi^{n+1}, \mu^{n+1}, \rho^{n+1}$ and $\boldsymbol{u}^{n+1}$ by solving

$$(3.42) \quad \frac{\phi^{n+1} - \phi^n}{\delta t} + \nabla \cdot (\boldsymbol{u}^{n+1}\phi^n) = M\Delta\mu^{n+1},$$

$$(3.43) \quad \mu^{n+1} = \lambda\left(-\epsilon\Delta\phi^{n+1} + \widehat{f}(\phi^n) + S_1(\phi^{n+1} - \phi^n)\right),$$

$$(3.44) \quad \rho^n\frac{\boldsymbol{u}^{n+1} - \boldsymbol{u}^n}{\delta t} - \nabla\cdot\nu^n D(\boldsymbol{u}^{n+1}) + \nabla(2p^n - p^{n-1}) + \rho^n(\boldsymbol{u}^n\cdot\nabla)\boldsymbol{u}^{n+1} + J^n\cdot\nabla\boldsymbol{u}^{n+1}$$
$$+ \phi^n\nabla\mu^{n+1} + \frac{1}{2}\frac{\rho^{n+1} - \rho^n}{\delta t}\boldsymbol{u}^{n+1} + \frac{1}{2}\nabla\cdot(\rho^n\boldsymbol{u}^n)\boldsymbol{u}^{n+1} + \frac{1}{2}\nabla\cdot J^n\boldsymbol{u}^{n+1} = 0,$$

with boundary conditions

$$(3.45) \quad \partial_{\boldsymbol{n}}\mu^{n+1} = 0, \quad \text{on } \Gamma,$$

$$(3.46) \quad \frac{\phi^{n+1} - \phi^n}{\delta t} + \boldsymbol{u}_\tau^{n+1}\cdot\nabla_\tau\phi^n = -\gamma\tilde{L}^{n+1}, \quad \text{on } \Gamma,$$

$$(3.47) \quad \boldsymbol{u}^{n+1}\cdot\boldsymbol{n} = 0, \quad \text{on } \Gamma,$$

$$(3.48) \quad \nu^n\partial_{\boldsymbol{n}}\boldsymbol{u}_\tau^{n+1} + \nu^n\ell(\phi^n)\boldsymbol{u}_s^{n+1} - \lambda\tilde{L}^{n+1}\nabla_\tau\phi^n = 0, \quad \text{on } \Gamma,$$

where

$$(3.49) \quad \tilde{L}^{n+1} := \epsilon\partial_{\boldsymbol{n}}\phi^{n+1} + g'(\phi^n) + S_2(\phi^{n+1} - \phi^n),$$

$$(3.50) \quad \rho^{n+1} = \frac{\rho_1 - \rho_2}{2}\widehat{\phi}^{n+1} + \frac{\rho_1 + \rho_2}{2},$$

$$(3.51) \quad \nu^{n+1} = \frac{\nu_1 - \nu_2}{2}\widehat{\phi}^{n+1} + \frac{\nu_1 + \nu_2}{2},$$

$$(3.52) \quad J^n = -M\frac{\rho_1 - \rho_2}{2}\nabla\mu^n.$$

*Step 2:* Update $p^{n+1}$ by solving

$$(3.53) \quad \Delta(p^{n+1} - p^n) = \frac{\chi}{\delta t}\nabla\cdot\boldsymbol{u}^{n+1},$$

with boundary conditions

(3.54) $$\partial_{\boldsymbol{n}} p^{n+1} = 0 \quad \text{on } \Gamma,$$

and $\chi = \frac{1}{2}\min(\rho_1, \rho_2)$.

**Theorem 3.** *Assuming $\boldsymbol{u}_w = 0$, and $S_1 \geq L_1/2$, $S_2 \geq L_2/2$, then the scheme (3.42)-(3.54) is energy stable in the sense that*

(3.55) $$E^{n+1} \leq E^n - \frac{\delta t}{2}\|\sqrt{\nu^n}D(\boldsymbol{u}^{n+1})\|^2 - \delta t M\|\nabla \mu^{n+1}\|^2 - \|\sqrt{\nu^n \ell(\phi^n)}\boldsymbol{u}_s^{n+1}\|_\Gamma^2 - \lambda \delta t \gamma \|\tilde{L}^{n+1}\|_\Gamma^2,$$

*where*

(3.56) $$E^n = \frac{1}{2}\|\sigma^n \boldsymbol{u}^n\|^2 + \lambda\left(\frac{\epsilon}{2}\|\nabla \phi^n\|^2 + (\widehat{F}(\phi^n), 1)\right) + \frac{\delta t^2}{2\chi}\|\nabla p^n\|^2 + \lambda (g(\phi^n), 1)_\Gamma.$$

*Proof.* (i) By taking the $L^2$ inner product of (3.44) with $2\delta t \boldsymbol{u}^{n+1}$, all the detailed derivations are same as Theorem 2 up to (3.35) except the boundary term as follows.

(3.57) $$\|\sigma^{n+1}\boldsymbol{u}^{n+1}\|^2 - \|\sigma^n \boldsymbol{u}^n\|^2 + \delta t\|\sqrt{\nu^n}D(\boldsymbol{u}^{n+1})\|^2 + \frac{\delta t^2}{\chi}(\|\nabla p^{n+1}\|^2 - \|\nabla p^n\|^2)$$
$$+ \frac{\delta t^2}{\chi}\|\nabla(p^n - p^{n-1})\|^2 \leq -2\delta t(\phi^n \nabla \mu^{n+1}, \boldsymbol{u}^{n+1}) + 2\delta t(\nu^n \partial_{\boldsymbol{n}} \boldsymbol{u}^{n+1}, \boldsymbol{u}^{n+1})_\Gamma.$$

For the boundary term in (3.57), using (3.48), we derive

(3.58) $$2\delta t(\nu^n \partial_{\boldsymbol{n}} \boldsymbol{u}^{n+1}, \boldsymbol{u}^{n+1})_\Gamma = 2\delta t(-\nu^n \ell(\phi^n)\boldsymbol{u}_s^{n+1} + \lambda \tilde{L}^{n+1}\nabla_\tau \phi^n, \boldsymbol{u}_\tau^{n+1})_\Gamma$$
$$= -2\delta t\|\sqrt{\nu^n \ell(\phi^n)}\boldsymbol{u}_s^{n+1}\|_\Gamma^2 + 2\delta t \lambda (\tilde{L}^{n+1}\nabla_\tau \phi^n, \boldsymbol{u}_\tau^{n+1})_\Gamma.$$

(ii) and (iii) By taking $L^2$ inner product of (3.42) with $2\delta t \mu^{n+1}$ and taking inner product of (3.43) with $-2(\phi^{n+1} - \phi^n)$, following same procedure as in proving Theorem 2, we get

(3.59) $$2(\phi^{n+1} - \phi^n, \mu^{n+1}) - 2\delta t\left(\boldsymbol{u}^{n+1}\phi^n, \nabla \mu^{n+1}\right) = -2\delta t M\|\nabla \mu^{n+1}\|^2.$$

and

(3.60) $$-2(\mu^{n+1}, \phi^{n+1} - \phi^n) = -\lambda \epsilon \left(\|\nabla \phi^{n+1}\|^2 - \|\nabla \phi\|^2 + \|\nabla(\phi^{n+1} - \phi^n)\|^2\right)$$
$$- 2\lambda(\widehat{f}(\phi^n) + S_1(\phi^{n+1} - \phi^n), \phi^{n+1} - \phi^n)$$
$$+ 2\lambda \epsilon(\partial_{\boldsymbol{n}}\phi^{n+1}, \phi^{n+1} - \phi^n)_\Gamma.$$

For the bounary term in (3.60), applying (3.46), we have

(3.61) $$2\lambda \epsilon(\partial_{\boldsymbol{n}}\phi^{n+1}, \phi^{n+1} - \phi^n)_\Gamma$$
$$= 2\lambda(\tilde{L}^{n+1}, \phi^{n+1} - \phi^n) - 2\lambda(g'(\phi^n) + S_2(\phi^{n+1} - \phi^n), \phi^{n+1} - \phi^n)_\Gamma$$
$$= 2\lambda \delta t(\tilde{L}^{n+1}, -\gamma \tilde{L}^{n+1} - \boldsymbol{u}_\tau^{n+1} \cdot \nabla_\tau \phi^n)_\Gamma - 2\lambda(g'(\phi^n) + S_2(\phi^{n+1} - \phi^n), \phi^{n+1} - \phi^n)_\Gamma$$
$$= -2\lambda \delta t \gamma \|\tilde{L}^{n+1}\|_\Gamma^2 - 2\lambda \delta t(\tilde{L}^{n+1}\nabla_\tau \phi^n, \boldsymbol{u}_\tau^{n+1})_\Gamma - 2\lambda(g'(\phi^n) + S_2(\phi^{n+1} - \phi^n), \phi^{n+1} - \phi^n)_\Gamma.$$

(iv) Summing up equations (3.57), (3.58), (3.59), (3.60) and (3.61), and applying the taylor expansion formulation (3.39) and (3.40), we get

$$\frac{1}{2}\|\sigma^{n+1}\boldsymbol{u}^{n+1}\|^2 - \frac{1}{2}\|\sigma^n \boldsymbol{u}^n\|^2 + \lambda\left(\frac{\epsilon}{2}\|\nabla\phi^{n+1}\|^2 - \frac{\epsilon}{2}\|\nabla\phi\|^2 + \frac{\epsilon}{2}\|\nabla(\phi^{n+1} - \phi^n)\|^2\right)$$
$$+ \lambda(\widehat{F}(\phi^{n+1}) - \widehat{F}(\phi^n), 1) + \lambda(g(\phi^{n+1}) - g(\phi^n), 1)_\Gamma$$
$$+ \frac{\delta t^2}{2\chi}\left(\|\nabla p^{n+1}\|^2 - \|\nabla p^n\|^2\right) + \frac{\delta t^2}{2\chi}\|\nabla p^{n+1} - \nabla p^n\|^2$$
$$\leq -\frac{\delta t}{2}\|\sqrt{\nu^n}D(\boldsymbol{u}^{n+1})\|^2 - \delta t M\|\nabla \mu^{n+1}\|^2 - \delta t\left\|\sqrt{\nu^n \ell(\phi^n)}\boldsymbol{u}_s^{n+1}\right\|_\Gamma^2 - \lambda\delta t\gamma\|\tilde{L}^{n+1}\|_\Gamma^2$$
$$- \lambda\left(S_1 - \frac{\widehat{f}'(\xi)}{2}, (\phi^{n+1} - \phi^n)^2\right) - \lambda\left(S_2 - \frac{g''(\zeta)}{2}, (\phi^{n+1} - \phi^n)^2\right)_\Gamma.$$

By the assumption $S_1 \geq L_1/2$, and $S_2 \geq L_2/2$, we get the desired energy estimate. □

**Remark 3.2.** *The equations (3.42)-(3.52) form a coupled system with nonlinearity through the term $\rho^{n+1}\boldsymbol{u}^{n+1}$ in (3.44). Note that for small Capillary numbers (corresponding to large $\lambda$, small $\nu$), the coupling can make the system stiff (cf. [2]). Nevertheless, the coupling system can be solved by either decoupling the system with a lagged velocity for the convective term in the phase equation or using a simple sub-iteration process. In our simulations, instead of using sub-iterations, we replace $\boldsymbol{u}^{n+1}$ by $\boldsymbol{u}^n$ in (3.42) and (3.46), then we can obtain $\phi^{n+1}$ and $\boldsymbol{u}^{n+1}$ by solving two decoupled linear elliptic equations. Furthermore, we choose small time steps to obtain the desired accuracy and avoid the spurious solutions. We call this scheme Linear Decoupled scheme with Explicit convection (LDE) scheme. Numerical results about the stability property of LDE scheme will be presented in Section 5.*

## 4. Spatial Discretization

In this section, we describe the spatial discretization using the Galerkin approach. As an example, we present a spectral Galerkin method for the CHNS coupled system in a 2-dimensional rectangular domain $\Omega = [0, L_x] \times [-1, 1]$. Finite element methods and spectral element methods for 2-D and 3-D domains can be built in a similar way.

**4.1. Weak formulations.** We take the LDS scheme for the CHNS system with GNBC and DCLC as an example to demonstrate the spatial discretization. In the first step, we solve (3.4)-(3.9), with (3.7) replaced by the discretized DCLC

$$(4.1) \qquad \frac{\phi^{n+1} - \phi^n}{\delta t} + \nabla \cdot (\boldsymbol{u}^n \phi^n) = -\gamma \tilde{L}^{n+1}.$$

Note that the discretized SCLC (3.7) is included as a special case ($1/\gamma = 0$) in the discretized DCLC (4.1). By moving the terms in (3.4)-(3.5) involves unknown variables to the left hand side and other terms to the right hand side, we get

$$(4.2) \qquad \frac{1}{\delta t}\phi^{n+1} - \delta t \nabla \cdot \left(\frac{(\phi^n)^2}{\rho^n}\nabla \mu^{n+1}\right) - M\Delta\mu^{n+1} = \frac{1}{\delta t}\phi^n - \nabla \cdot (\boldsymbol{u}^n \phi^n),$$

$$(4.3) \qquad -\epsilon\Delta\phi^{n+1} + S_1 \phi^{n+1} - \frac{1}{\lambda}\mu^{n+1} = S_1 \phi^n - \widehat{f}(\phi^n),$$

The boundary conditions are

(4.4) $$\partial_{\boldsymbol{n}}\mu^{n+1} = 0, \quad \text{on } \Gamma,$$

(4.5) $$\left(\frac{1}{\gamma\delta t} + S_2\right)\phi^{n+1} + \epsilon\partial_{\boldsymbol{n}}\phi^{n+1} = \left(\frac{1}{\gamma\delta t} + S_2\right)\phi^n - \frac{1}{\gamma}\nabla\cdot(\boldsymbol{u}^n\phi^n) - g'(\phi^n), \quad \text{on } \Gamma.$$

The corresponding weak formulation for the above equations reads as follows:

Find $\mu^{n+1} \in H^1(\Omega), \phi^{n+1} \in H^1(\Omega)$, such that for any $\omega \in H^1(\Omega)$ and $\varphi \in H^1(\Omega)$,

(4.6) $$\epsilon(\nabla\phi^{n+1}, \nabla\varphi) + S_1(\phi^{n+1}, \varphi) + c_s \int_\Gamma \phi^{n+1}\varphi - \frac{1}{\lambda}(\mu^{n+1}, \varphi) = (S_1\phi^n - \widehat{f}(\phi^n), \varphi) + \int_\Gamma \phi_\Gamma^n \varphi,$$

and

(4.7) $$\frac{1}{\delta t}(\phi^{n+1}, \omega) + \delta t\left(\frac{(\phi^n)^2}{\rho^n}\nabla\mu^{n+1}, \nabla\omega\right) + M(\nabla\mu^{n+1}, \nabla\omega) = (\phi_r^n, \omega).$$

Here, $c_s = \frac{1}{\delta t\gamma} + S_2$, $\phi_r^n = \frac{1}{\delta t}\phi^n - \nabla\cdot(\boldsymbol{u}^n\phi^n)$, $\phi_\Gamma^n = \frac{1}{\gamma}\phi_r^n - \widehat{g}'(\phi^n) + S_2\phi^n$. Note that when $\gamma \to \infty$, the above formulation handles the system with static contact line condition.

In the second step, we solve equations (3.12)-(3.14) with (3.14) is replaced by

(4.8) $$\nu^n\partial_{\boldsymbol{n}}\boldsymbol{u}_\tau^{n+1} + \nu^n\ell(\phi^n)\boldsymbol{u}_s^{n+1} - \lambda\tilde{L}^{n+1}\nabla_\tau\phi^n = 0, \quad \text{on } \Gamma.$$

The corresponding weak formulation is:

Find $\boldsymbol{u}^{n+1} \in V_{\boldsymbol{u}} := H^1(\Omega) \times H_0^1(\Omega)$, such that for any $\boldsymbol{v} \in V_{\boldsymbol{u}}$,

(4.9) $$\begin{aligned}\frac{1}{\delta t}(\rho^n\boldsymbol{u}^{n+1}, \boldsymbol{v}) + ((\rho^n\boldsymbol{u}^n + J^n)\cdot\nabla\boldsymbol{u}^{n+1}, \boldsymbol{v}) + \frac{1}{2}(\nu^n D(\boldsymbol{u}^{n+1}), D(\boldsymbol{v})) + \int_\Gamma \nu^n\ell(\phi^n)\boldsymbol{u}^{n+1}\cdot\boldsymbol{v}_\tau \\ + \frac{1}{2}\left(\left(\frac{\rho^{n+1} - \rho^n}{\delta t} + \nabla\cdot(\rho^n\boldsymbol{u}^n) + \nabla\cdot J^n\right)\boldsymbol{u}^{n+1}, \boldsymbol{v}\right) \\ = (\boldsymbol{u}_r^n, \boldsymbol{v}) + \int_\Gamma \boldsymbol{u}_\Gamma^n\cdot\boldsymbol{v}_\tau ds,\end{aligned}$$

where $\boldsymbol{u}_r^n = \frac{1}{\delta t}\rho^n\boldsymbol{u}^n - \phi^n\nabla\mu^{n+1} - \nabla(2p^n - p^{n-1})$, $\boldsymbol{u}_\Gamma^n = \lambda\tilde{L}^{n+1}\nabla\phi^n + \nu^n\ell(\phi^n)\boldsymbol{u}_w$.

For the last step, we solve (3.16) for $p^{n+1}$. The corresponding weak form is:

Find $p^{n+1} \in H_c^1(\Omega) := \{p : p \in H^1(\Omega), \int_\Omega p\,dx = 0\}$, such that for any $q \in H_c^1(\Omega)$

(4.10) $$(\nabla p^{n+1}, \nabla q) = (\nabla p^n, \nabla q) - \frac{\chi}{\delta t}(\nabla\cdot\boldsymbol{u}^{n+1}, q).$$

The weak formulation for the LDE scheme is similar, except that the term $\delta t\nabla\cdot\left(\frac{(\phi^n)^2}{\rho^n}\nabla\mu^{n+1}\right)$ in (4.2) does not show up in the weak form of LDE scheme.

4.2. **Spectral Galerkin approximation and solution procedure.** We assume the system in $x$ direction is periodic, only the top and bottom boundaries take the GNBC and DCLC. Similarly as in [54], we use real Fourier bases in $x$ direction and Legendre quasi-orthogonal polynomials in $y$ direction. More precisely, we use

(4.11) $$F_m := \text{span}\{E_k(x), 0 \leq k < m\}, \quad P_k = \text{span}\{\varphi_j(y) : 0 \leq j < k\}$$

as the basis set for the $x$-direction and $y$-direction correspondingly, where

$$(4.12) \qquad E_k(x) = \begin{cases} \cos(k\pi x/L_x), & k \text{ even}, \\ \sin((k+1)\pi x/L_x), & k \text{ odd}, \end{cases}$$

$$(4.13) \qquad \varphi_0(y) = \frac{1+x}{2}, \quad \varphi_1(y) = \frac{1-x}{2}, \quad \varphi_j(y) = L_j(y) - L_{j-2}(y), \quad j = 2, 3, \ldots,$$

and $L_j(y)$ is the Legendre polynomial of degree $j$. This is a direct extension of nearly orthogonal bases proposed by Shen [45]. For given $n_x, n_y$, we take $F_{n_x} \otimes P_{n_y}$ as the approximation space for $\mu^{n+1}$ and $\phi^{n+1}$. For the Navier–Stokes equation, the velocity in $x$-component satisfies the GNBC, which is a Robin type boundary condition, while the component in $y$-direction satisfies Dirichlet boundary condition. The Robin type boundary condition is treated naturally in the weak form, the Dirichlet boundary condition is imposed on the approximation space. So the Galerkin approximation space for $V_{\boldsymbol{u}}$ is $V^{\boldsymbol{u}}_{n_x,n_y} := \{ (u,v) \,|\, u \in F_{n_x} \otimes P_{n_y}, v \in F_{n_x} \otimes P^0_{n_y} \}$, where $P^0_{n_y} = \text{span}\{\varphi_j(y), j = 2, .., n_y\}$. The approximation space for pressure is $V^p_{n_x,n_y} = F_{n_x} \otimes P_{n_y} \backslash C := \text{span}\left\{E_l(x)\varphi_j(y) : 0 \leq l < n_x, 0 \leq j < n_y, l^2 + j^2 \neq 0\right\}$. One of the advantages of this basis set is that equations with constant coefficients all lead to sparse linear algebra systems, which can be solved with optimal complexity. For variable-coefficient equations, we use a preconditioned CG method (BiCGSTAB) with matrices for corresponding constant-coefficient systems as preconditioners. All the nonlinear terms are first evaluated in physical space and then convert into spectral space using fast spectral transform. No anti-aliasing rule is used for the nonlinear terms, since we found when spatial resolution is good enough, anti-aliasing has almost no effects on the numerical results.

## 5. Numerical Simulations

In this section, we present some numerical results to validate our proposed schemes. We first perform a convergence test, and then present numerical results of several typical cases. There are a lot of parameters in the system. If not explicit specified, the model parameters take default values given below:

$$(5.1) \quad \begin{aligned} & L_x = 6, \quad \lambda = 1.2, \quad M = 0.01, \quad \gamma = 100, \quad \epsilon = 0.05, \quad \theta_s = \pi/3, \\ & \rho_1 = 1, \quad \rho_2 = 0.9, \quad \nu_1 = 1, \quad \nu_2 = 1.1, \quad \nu l(\phi) = 1/0.19, \quad u_w^\pm = \pm 0.2, \end{aligned}$$

where $u_w^\pm$ is the velocities of uppper and bottom plates.

5.1. **Convergence and performance test.** We first present the results of convergence test. Here we let the initial velocity field take the profile of Couette flow, and set the initial value of $\phi$ as

$$(5.2) \qquad \phi_0(x,y) = \tanh\left(\frac{1}{\sqrt{2}\epsilon}\left(\frac{L_x}{4} - \left|x - \frac{L_x}{2}\right|\right)\right).$$

First, we use two schemes proposed in this paper to solve the system, with 129 Fourier modes, and 48 Legendre modes and $\delta t = 0.01$. The numerical results of the two schemes are similar, and in both schemes, the volume of $\phi$ are kept up to machine accuracy. In Fig. 5.1, we show the numerical solutions of LDE scheme at $T = 2$.

To test the convergence for time step $\delta t$, we use a spatial resolution $n_x = 257, n_y = 64$ for a series of temporal resolution $\delta t = 0.004, 0.002, 0.001, 0.0005, 0.00025$. Since the exact solution is unknown, a solution calculated with a very small time step-size $\delta t = 1.25 \times 10^{-4}$ is used as the exact solution. In Table 5.1, we present the $L^2$ error of the velocity and phase variable of the two numerical schemes. We see from this table that both schemes can achieve first order accuracy in time.

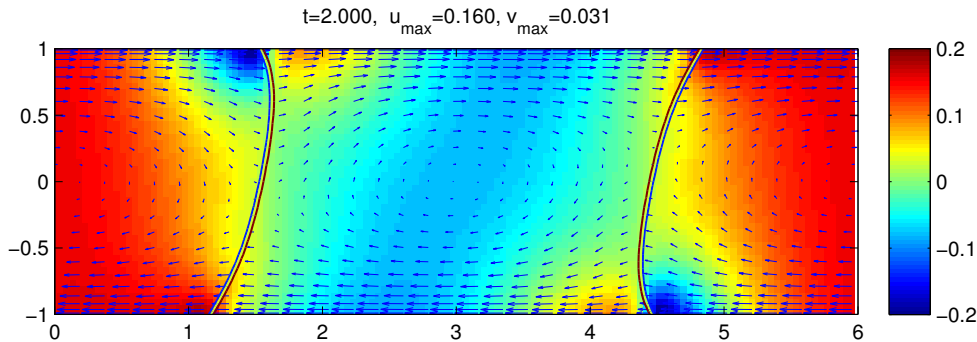

FIGURE 5.1. The contour line of $\phi(x,t) = \{-0.2, 0, 0.2\}$, quiver plot of velocity $u, v$ (arrows) and pressure fields (background color) at $t = 2$ of the numerical results of LDE scheme using 129 Fourier modes and 48 Legendre modes, $\delta t = 0.01$. Model parameters are given in (5.1).

| $\delta t$ | $Error_u^{\text{LDS}}$ | Order | $Error_\phi^{\text{LDS}}$ | Order | $Error_u^{\text{LDE}}$ | Order | $Error_\phi^{\text{LDE}}$ | Order |
|---|---|---|---|---|---|---|---|---|
| 0.004 | 3.424(-2) | | 2.834(-2) | | 3.398(-2) | | 3.851(-2) | |
| 0.002 | 1.778(-2) | 0.946 | 1.396(-2) | 1.021 | 1.767(-2) | 0.943 | 1.994(-2) | 0.950 |
| 0.001 | 8.602(-3) | 1.047 | 6.570(-3) | 1.088 | 8.567(-3) | 1.045 | 9.631(-3) | 1.050 |
| 0.0005 | 3.755(-3) | 1.196 | 2.830(-3) | 1.215 | 3.744(-3) | 1.194 | 4.199(-3) | 1.197 |
| 0.00025 | 1.263(-3) | 1.572 | 9.447(-4) | 1.583 | 1.260(-3) | 1.571 | 1.412(-3) | 1.573 |

TABLE 5.1. The numerical errors for velocity $\boldsymbol{u}$ and phase variable $\phi$ at $t = 0.2$ for both of LDS and LDE schemes using different temporal resolutions. The numerical results using $\delta t = 0.000125$ is taken as the exact solution to calculate the $L^2$ errors. Model parameters given by (5.1) are used where $\rho_1 = 1, \rho_2 = 0.9$.

We further test the temporal convergence of LDE and LDS schemes for a large density/viscosity ratio. All parameters are from (5.1) except $\rho_1 = 0.1, \rho_2 = 10, \nu_1 = 0.1, \nu_2 = 10$. The numerical errors are shown in Table 5.2, where the first order temporal convergence is achieved as well.

To test the convergence for space discretization, we use $n_y = 64$, $\delta t = 0.0001$ for a series of $n_x$ starting from 113 with the increment 16. The $L^2$ error is calculated at time $t = 2$, with a reference solution obtained using resolution $n_x = 257, n_y = 64, \delta t = 0.0001$. Note that for $t = 2$ the system is very close to the steady state, so the error related to temporal discretization is negligible. To test the convergence in $y$, we use $n_x = 257$, $\delta t = 0.0001$ for a series of $n_y$ starting from 16 with an increment 4. The $L^2$ error is again calculated at $t = 2$. The reference solution is still the obtained numerical solution using resolutions of $n_x = 257, n_y = 52, \delta t = 0.0001$. Since the convergence behavior of the LDE and LDS schemes are almost identical, we show only the results of LDS scheme in Fig. 5.2, from which, we see that the proposed numerical schemes can achieve spectral accuracy in $L^2$ norm.

| $\delta t$ | $Error_u^{\text{LDS}}$ | Order | $Error_\phi^{\text{LDS}}$ | Order | $Error_u^{\text{LDE}}$ | Order | $Error_\phi^{\text{LDE}}$ | Order |
|---|---|---|---|---|---|---|---|---|
| 0.008 | 1.847(-2) | | 5.398(-2) | | 1.922(-2) | | 2.581(-2) | |
| 0.004 | 1.111(-2) | 0.733 | 2.813(-2) | 0.940 | 1.142(-2) | 0.751 | 1.303(-2) | 0.986 |
| 0.002 | 5.324(-3) | 1.062 | 1.448(-2) | 0.958 | 5.377(-3) | 1.087 | 6.603(-3) | 0.981 |
| 0.001 | 2.318(-3) | 1.200 | 5.300(-3) | 1.450 | 2.313(-3) | 1.217 | 2.944(-3) | 1.165 |
| 0.0005 | 9.174(-4) | 1.337 | 1.746(-3) | 1.602 | 7.764(-4) | 1.575 | 1.002(-3) | 1.555 |

TABLE 5.2. The numerical errors for velocity $\boldsymbol{u}$ and phase variable $\phi$ at $t = 0.2$ for both of LDS and LDE schemes using different temporal resolutions for a large density/viscosity ratio case. The numerical results using $\delta t = 0.00025$ is taken as the exact solution to calculate the $L^2$ errors. The parameters given by (5.1) are used, except that here we take a large density and viscosity ratio: $\rho_1 = 0.1, \rho_2 = 10, \nu_1 = 0.1, \nu_2 = 10$.

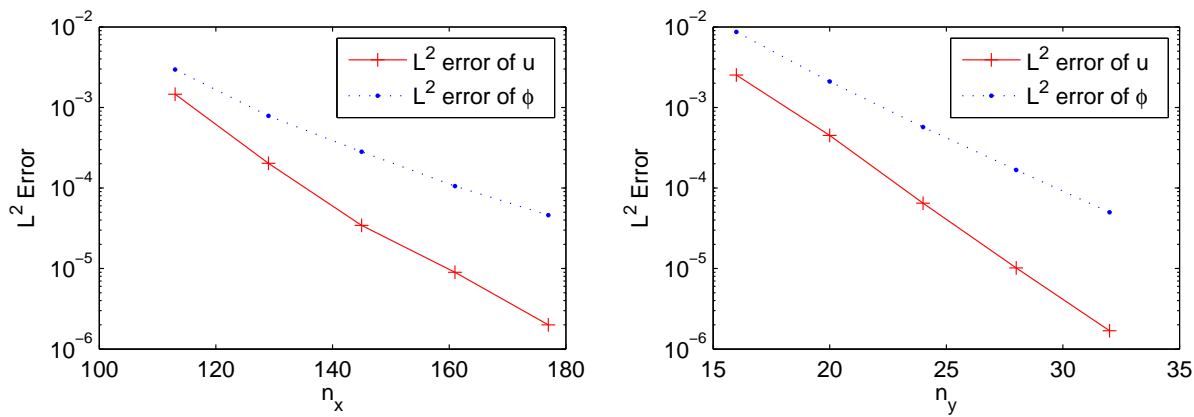

FIGURE 5.2. The convergence of spatial discretization of LDS scheme. The results on a fine grid $n_x = 257, n_y = 64$ is taken as the reference solution for the left part. Results on a fine grid $n_x = 257, n_y = 52$ is taken as the reference solution for the right part. Model parameters given by (5.1) are used. Spectral accuracy is achieved.

Next, we show the evolutions of the total free energy in Fig. 5.3 with respect to three time step sizes $\delta t = 0.04, 0.01, 0.0025$ for both LDS and LDE schemes. Here, we let $u_w^\pm = 0$, i.e. no input energy from outside. All other model parameters are given in (5.1). In Fig. 5.3 (a), we observe that all six energy curves decay, which confirm that our algorithms are energy stable. Meanwhile, the energy curve with larger time steps show considerable deviation away from the energy curve with smaller time steps. That means the induced numerical errors with larger time steps are higher. In Fig. 5.3 (b), we present the evolution of kinetic energy $E_k$ for LDS scheme with three time steps.

We observe that the flow kinetic energy first increase quickly to a peak value then gradually decay to zero.

The convergence with respect to the thickness paramter $\epsilon$ is one of the key issues of phase field model. We numerically verify this property for a case where $L_x = 4, u_w^\pm = 0, M = 1, \rho_1 = 1, \rho_2 = 0.8, \nu_1 = 1, \nu_2 = 2.0, \theta_s = \pi/6$, and other parameters take values given by (5.1). The inital profile of the phase field variable is a drop (half disk) with radius 1 setting on the bottom plate. We take a series of values of $\epsilon$. For each one, we run the LDE scheme enough long time to get close to the steady state. We then measure the contact angle by looking at $x$-coordinates of the contact points and radius of the drop. The result is given in Fig. 5.4. If use the results of $\epsilon \leq 0.06$ to fit the convergence order, we get the convergence of order $O(\epsilon^r)$ with $r \approx 1.6$, but if we use $\epsilon$ much closer to 0, we get a order close to 2. For example, if we use the results of $\epsilon \leq 0.015$ to fit convergence order, the result is about 1.9.

At the end of this subsection, we present the stability property of using LDE scheme to decouple the nonlinear scheme presented in Section 3.2. In Table 5.3, we show the maximum stable time step size for LDE scheme with spatial resolution $n_x = 65, n_y = 24$ to solve systems with different values of $\gamma, \lambda$ and other model parameters are given in (5.1). We observe that, The parameter $\gamma$ and spatial resolution have only very small effects on the maximum allowed time steps. But $\lambda$ has a strong effect on the stability of the LDE scheme. For larger value of $\lambda$ the maximum stable time step is much smaller.

|  | $n_x \times n_y = 95 \times 32$ | | | | $n_x \times n_y = 127 \times 48$ | | | |
| --- | --- | --- | --- | --- | --- | --- | --- | --- |
| $\lambda$ | 1 | 4 | 16 | 64 | 1 | 4 | 16 | 64 |
| $\gamma = 1$ | 8.1(-2) | 2.2(-2) | 5.2(-3) | 7.5(-4) | 8.0(-2) | 2.1(-2) | 5.2(-3) | 7.5(-4) |
| $\gamma = 4$ | 9.0(-2) | 2.4(-2) | 5.2(-3) | 7.5(-4) | 9.1(-2) | 2.3(-2) | 5.3(-3) | 7.5(-4) |
| $\gamma = 16$ | 9.1(-2) | 2.6(-2) | 5.6(-3) | 7.5(-4) | 9.2(-2) | 2.6(-2) | 5.4(-3) | 7.5(-4) |
| $\gamma = 64$ | 9.2(-2) | 2.6(-2) | 5.9(-3) | 7.7(-4) | 9.2(-2) | 2.6(-2) | 5.5(-3) | 7.7(-4) |

TABLE 5.3. Maximum time step size allowed for different values of $\gamma, \lambda$ and spatial resolution for the LDE scheme. Model parmaters are given in (5.1).

5.2. **A oil-water mix in shear flow without gravity.** In this example, we simulate the dynamics of a shear drop with small density and viscosity difference. The first phase (ambient fluid) has a slightly larger density and smaller viscosity, may be regarded as water, while the second phase (drop) has a slightly smaller density and larger viscosity, which may be regarded as oil. The model parameters are $u_w^\pm = \pm 2$, $\theta_s = \pi/6$, $\rho_1 = 1, \rho_2 = 0.8$ $\nu_1 = 1, \nu_2 = 2$, $\lambda = 12$, $M = 0.001$, and discretization parameters are $n_x = 197$, $n_y = 48$, $\delta t = 0.005$.

We plot the contour line of the phase-field variable $\{x : \phi(x,t) = 0\}$ together with the velocity field in Fig. 5.5. Initially, the oil drop has a smaller volume, and set in the middle of the bottom plate, as shown in the first figure of Fig. 5.5. The numerical results show that, as the bottom plate moves, the oil phase will also move with it but at a much smaller velocity. This means it slips on the bottom plate. As time goes on, the contact region of the oil phase with the bottom boundary gets smaller and smaller. The oil drop eventually gets off the bottom boundary at $t \approx 4.5$, and moves toward the center of the channel. We also simulate a larger static contact angle $\theta_s = \frac{2\pi}{3}$ case, in which the oil phase is harder to get off the bottom boundary, which is shown in Fig. 5.6.

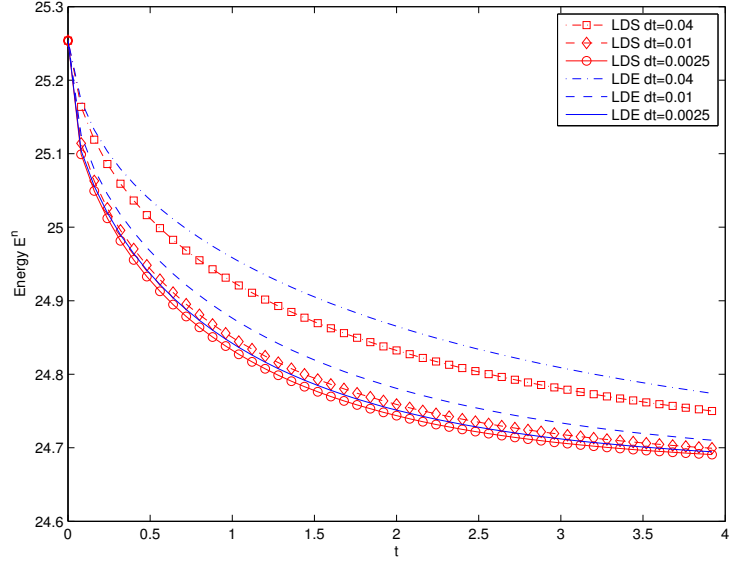

(a)

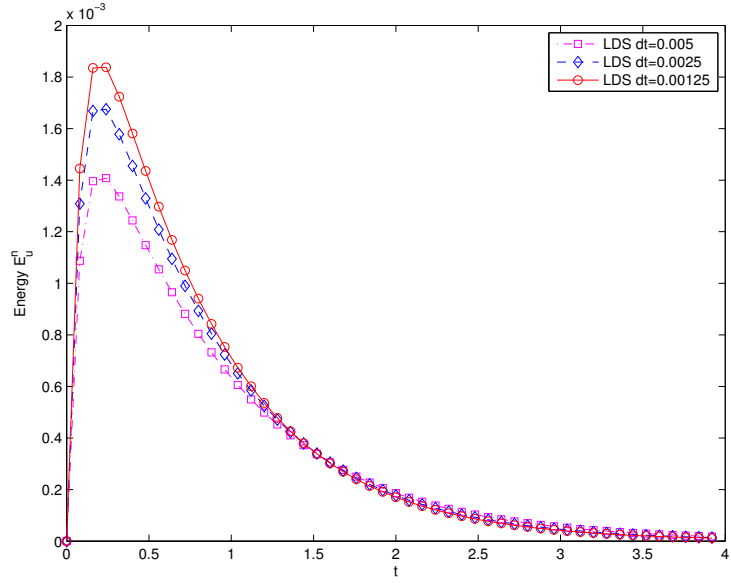

(b)

FIGURE 5.3. The energy dissipation properties of the LDS and LDE scheme. (a): energy dissipation using LDE and LDS scheme with different time stepsizes; (b) The evolution of kinetic energy $\|\sqrt{\rho}\boldsymbol{u}\|^2$ using LDS scheme with different time stepsizes. The spatial resolution is $n_x = 255, n_y = 64$, $u_w^\pm = 0$, other model parameters are given by (5.1).

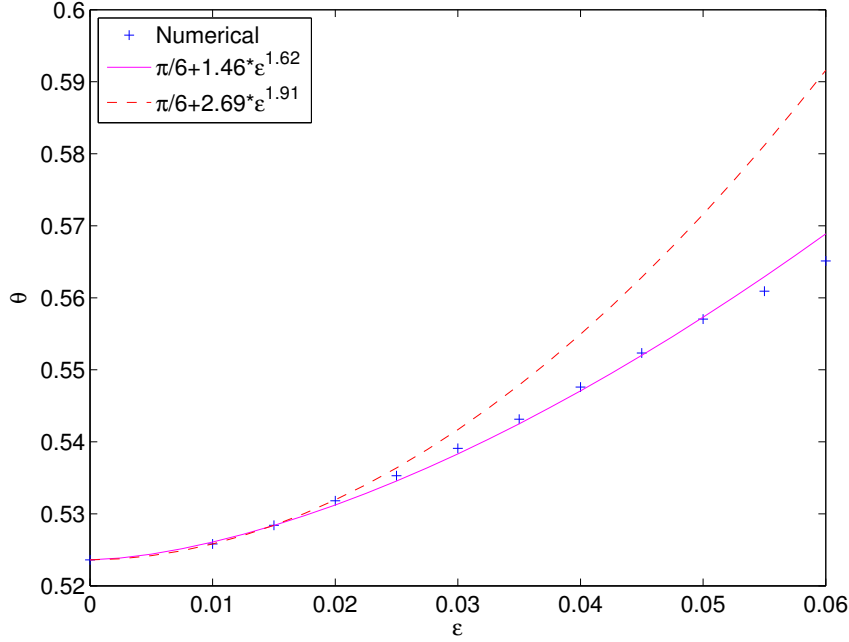

FIGURE 5.4. The convergence of CHNS model with respect to $\epsilon$. In this simulation, we take $L_x = 4, u_w^\pm = 0, M = 1, \rho_1 = 1, \rho_2 = 0.8, \nu_1 = 1, \nu_2 = 2.0$, other model parameters are given by (5.1).

**5.3. A droplet slides down a 45-degree slope.** In this case, we present the dynamical behaviors of a water drop surrounded by air. The drop will slide down along a 45-degree slope. We still use the rectangular domain, but a gravitational force $\rho \mathbf{g} = \rho(g_x, g_y)$ is added to the right hand side of the momentum equations of (2.7). We use following parameters

$$(5.3) \quad u_w^\pm = 0, \quad \mathbf{g} = (10, -10), \quad \theta_s = \frac{\pi}{6}, \quad \rho_1 = 10^{-3}, \quad \rho_2 = 1, \quad \nu_1 = 10^{-2}, \nu_2 = 1, \quad \gamma = 500.$$

We solve the system using the LDS scheme with grid points $n_x = 127, n_y = 64$, and $\delta t = 0.005$. Fig. 5.7 illustrates the dynamics of the drop until $t = 3.5$. Initially, at $t = 0$, the drop (half disk) is set at the bottom plate with the contact angle $\pi/2$ and all other fields are at rest. We observe the contact angle is already adjusted close to the prescribed angle $\theta_s$ at $t = 0.5$. Due to the gravity force, as time evolves, the bubble will slip along the bottom plate while keeping the contact angle almost unchanged.

**5.4. Air bubble rising in water.** As the last example, we simulate the air bubble rising in water. We use following model parameters

$$(5.4) \quad \begin{aligned} &L_x = 1, \ \epsilon = 0.02, \ u_w^\pm = 0, \ \mathbf{g} = (0, -100), \rho_1 = 1, \rho_2 = 0.002, \\ &\nu_1 = 1, \nu_2 = 0.01, \gamma = 10^6, M = 10^{-4}. \end{aligned}$$

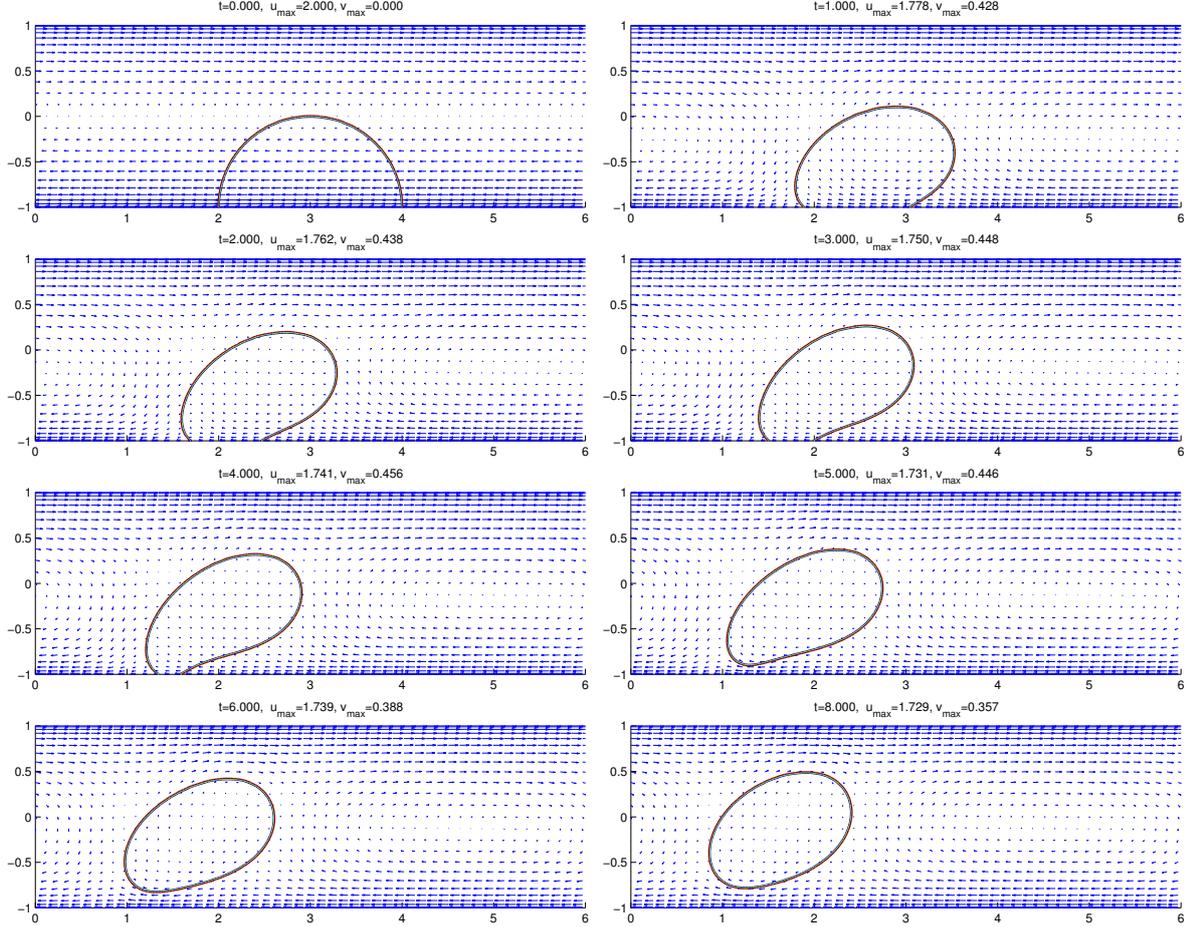

FIGURE 5.5. Oil-water mixture in shear flow without gravity. Results of LDE scheme with $n_x = 197$, $n_y = 48$, $\delta t = 0.005$. We take $u_w^\pm = \pm 2$, $\theta_s = \pi/6$, $\rho_2 = 0.8$, $\nu_2 = 2$, $\lambda = 12$, $M = 0.001$, and other equation parameters are given in (5.1).

We solve the system by using the LDS scheme with grid parameters $n_x = 129, n_y = 128$, $\delta t = 0.0001$. We compare the drop dynamics using two different contact angles of $\theta_s = \pi/6$ and $3\pi/4$ in Fig. 5.8 and Fig. 5.9, respectively. We notice that with an acute static contact angle, the air bubble is relatively easier to get off the bottom, shown in Fig. 5.8. With a obtuse static contact angle, the air bubble is harder to get off the bottom as a whole, a small part will be retained by the bottom plate, shown in Fig. 5.9. For both cases, the air bubble keeps a convex near-ball shape, since in the parameter setting (5.4), the corresponding Reynolds number is very small ($Re_1 = \rho_1 L |\boldsymbol{u}_1|/\nu_1 \approx 1, Re_2 \approx 0.2$). For larger Reynolds number, the air bubble keeps a non-convex shape during the rising process, but the process of getting off from the bottom plate is similar to the low Reynolds number case, see Fig. 5.10 for such a case.

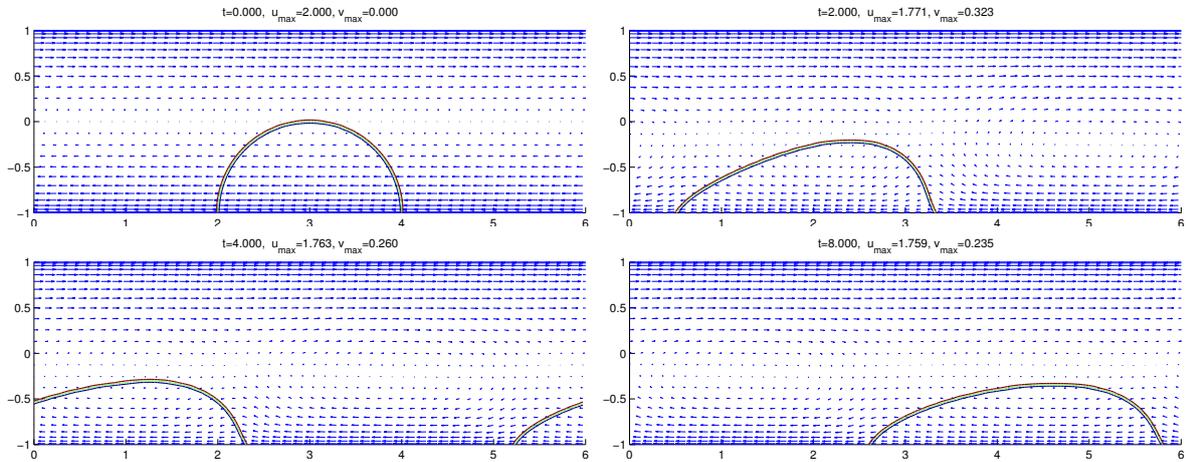

FIGURE 5.6. Oil-water mixture in shear flow without gravity. We take $u_w^\pm = \pm 2$, $\theta_s = 2\pi/3$, $\rho_2 = 0.8$, $\nu_2 = 2$, $\lambda = 12$, $M = 0.001$, and other equation parameters are given in (5.1). The results are of LDE scheme using $n_x = 197$, $n_y = 48$, $\delta t = 0.005$.

## 6. Summary

In this paper, we couple the variable density/viscosity to the classical phase-field model for the moving contact line problem. The model consists of incompressible Navier–Stokes equations with a generalized Navier boundary condition and the Cahn–Hilliard equation with the dynamic contact line condition, and satisfies an energy dissipation law. We then construct two energy stable, efficient schemes to solve the proposed nonlinear, coupled system. Various numerical simulations are presented to verify the accuracy and efficiency of the proposed model and numerical schemes.

## Acknowledgments

The work of H. Yu is partially supported by National Key Basic Research Program of China under Grants 2015CB856003 and the National NSF of China under Grants 91530322, 11101413, 11371358. X. Yang is partially supported by NSF DMS-1200487 and NSF DMS-1418898. X. Yang thanks the hospitality of Institute of Computational Mathematics of Chinese Academy of Science during his visit.

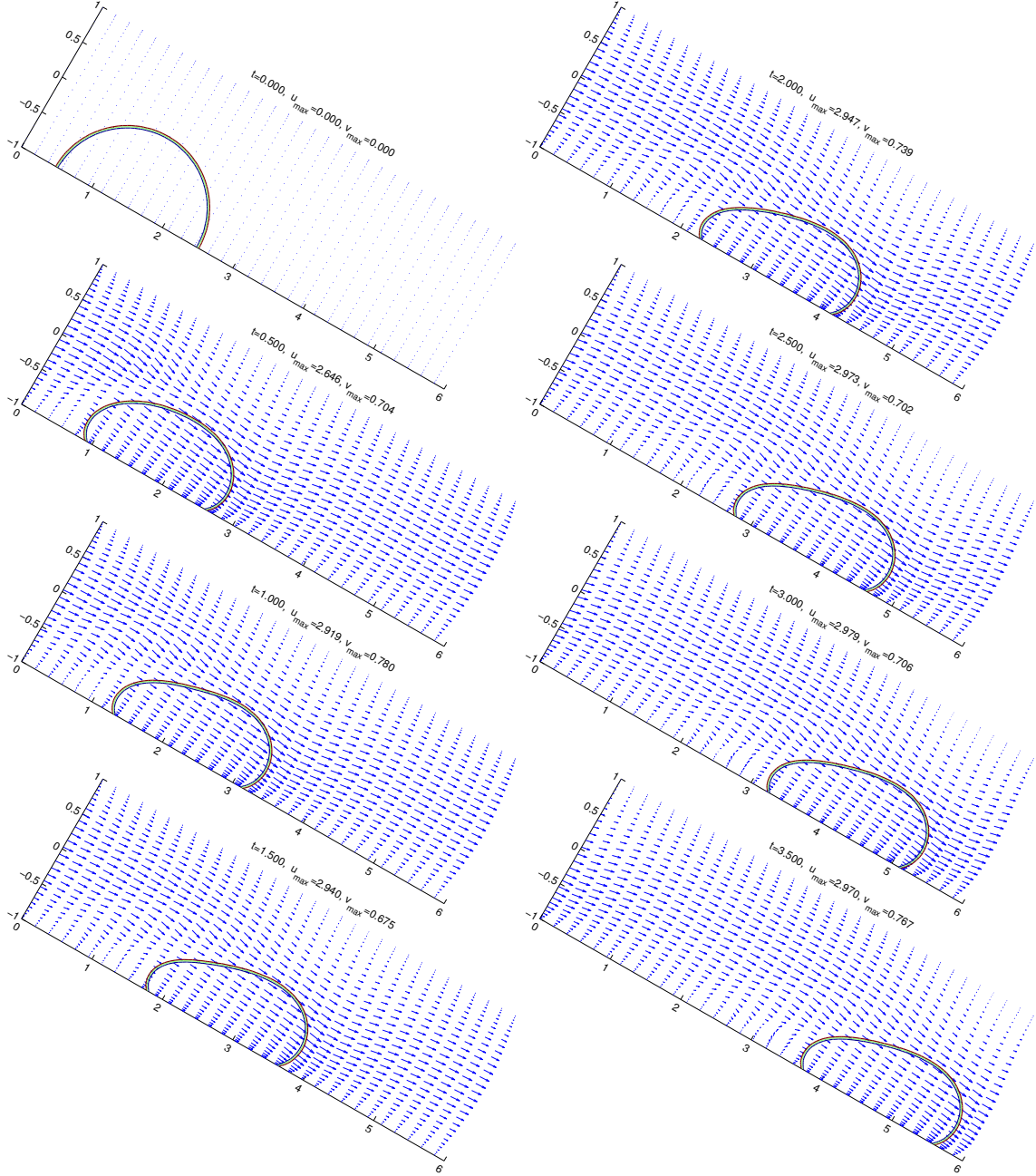

FIGURE 5.7. A water droplet slides down a 45-degree slope. The model parameters are given in (5.3), the parameters not specified in (5.3) take default values given in (5.1). The LDS scheme with $n_x = 127, n_y = 64$, $\delta t = 0.005$ is used to calculate the results.

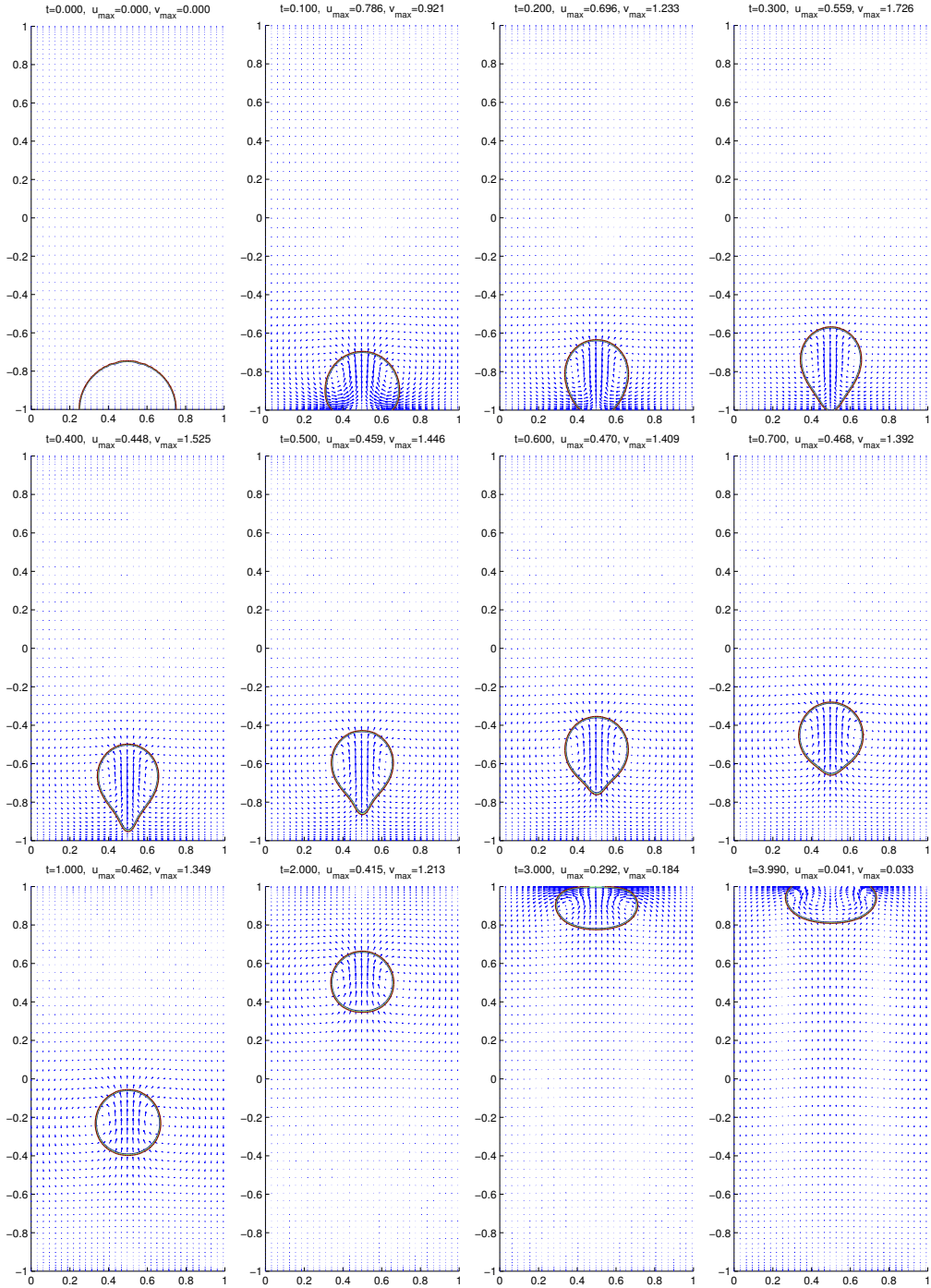

FIGURE 5.8. An air bubble with static contact angle $\theta_s = \pi/6$ rising. Other equation parameters are given in (5.4), or (5.1) if not specified in (5.4).

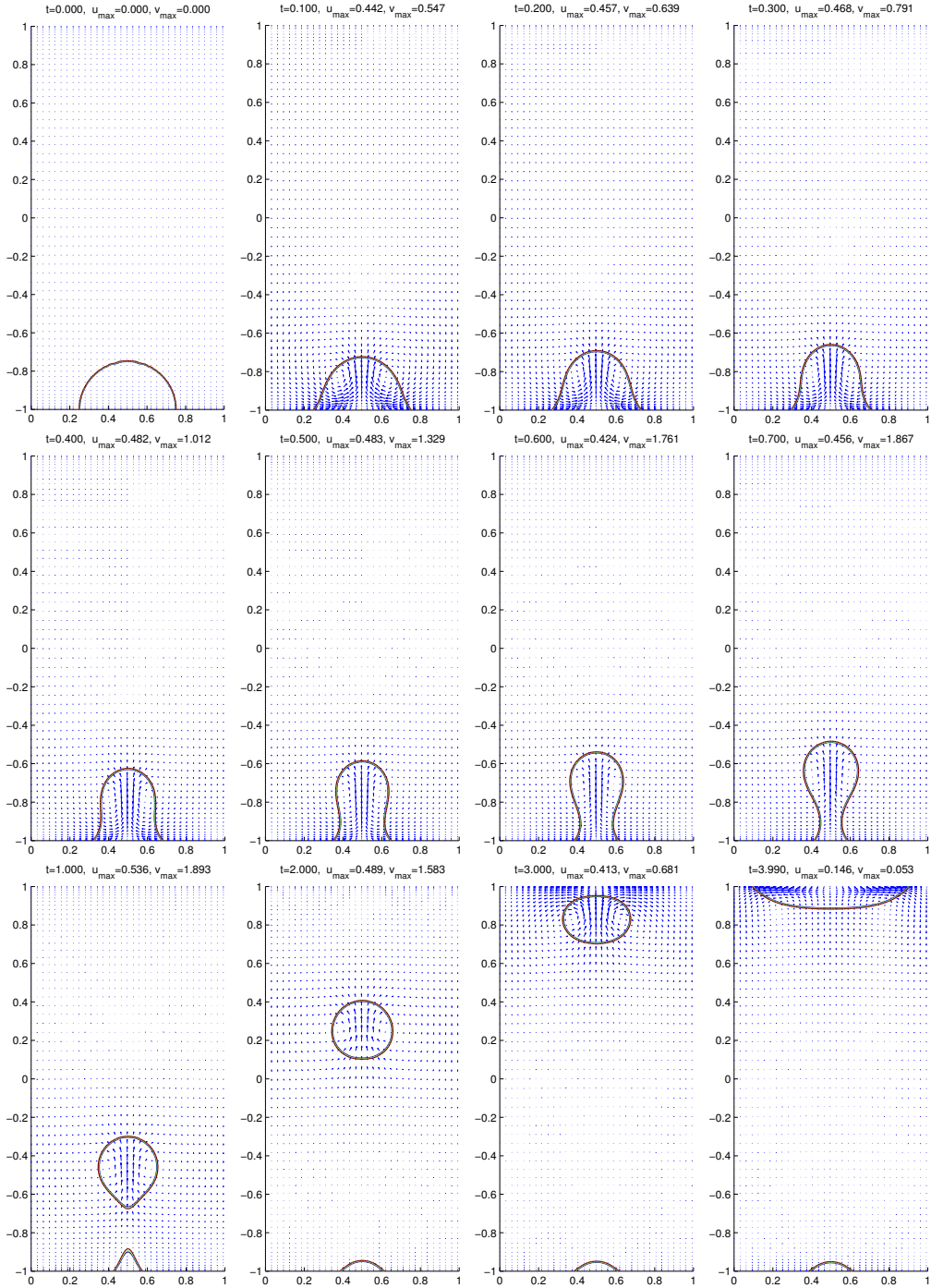

FIGURE 5.9. An air bubble with static contact angle $\theta_s = 3\pi/4$ rising. Other equation parameters are given in (5.4), or (5.1) if not specified in (5.4).

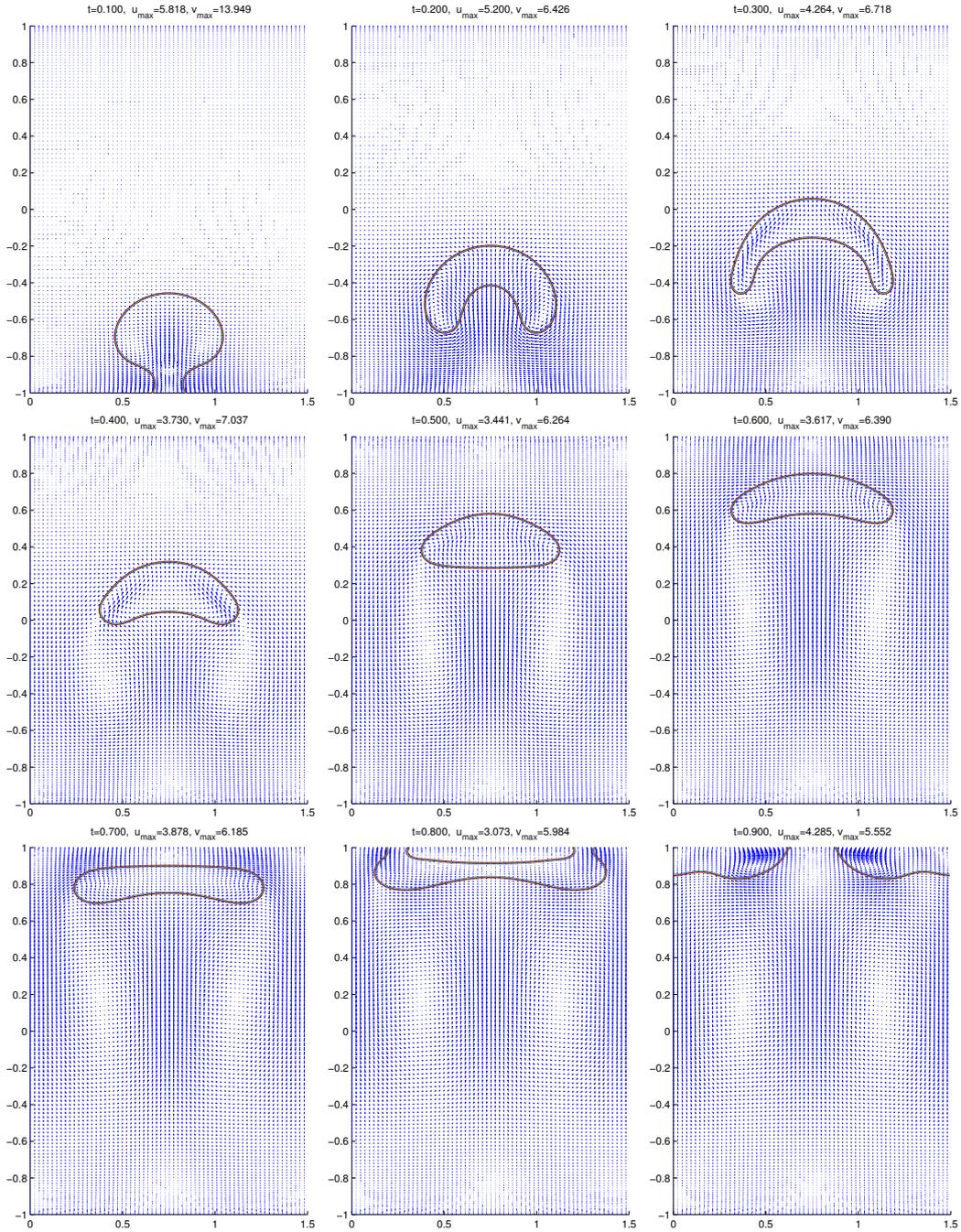

FIGURE 5.10. An air bubble rising with a non-convex shape. $L_x = 1.5, \epsilon = 0.02, u_w^\pm = 0, \theta_s = 2\pi/3, \mathbf{g} = (0, -100), \rho_1 = 1, \rho_2 = 0.002, \nu_1 = 0.032, \nu_2 = 0.00075, \gamma = 10^6, \lambda = 0.4$.